\documentclass[twocolumn,english]{IEEEtran}
\usepackage[T1]{fontenc}
\usepackage{babel}
\usepackage{float}
\usepackage{amsthm}
\usepackage{amsmath}
\usepackage{graphicx}
\usepackage[unicode=true,
 bookmarks=true,bookmarksnumbered=true,bookmarksopen=true,bookmarksopenlevel=1,
 breaklinks=false,pdfborder={0 0 0},backref=false,colorlinks=false]
 {hyperref}
\hypersetup{pdftitle={Your Title},
 pdfauthor={Your Name},
 pdfpagelayout=OneColumn, pdfnewwindow=true, pdfstartview=XYZ, plainpages=false}
\usepackage{breakurl}

\makeatletter

\floatstyle{ruled}
\newfloat{algorithm}{tbp}{loa}
\providecommand{\algorithmname}{Algorithm}
\floatname{algorithm}{\protect\algorithmname}

 \usepackage{algolyx}
\theoremstyle{plain}
\newtheorem{thm}{\protect\theoremname}
\theoremstyle{plain}
\newtheorem{lem}[thm]{\protect\lemmaname}
\theoremstyle{remark}
\newtheorem{rem}[thm]{\protect\remarkname}
\theoremstyle{plain}
\newtheorem{prop}[thm]{\protect\propositionname}

\ifCLASSOPTIONcompsoc
\usepackage[caption=false,font=normalsize,labelfont=sf,textfont=sf]{subfig}
\else
\usepackage[caption=false,font=footnotesize]{subfig}
\fi

\usepackage{cite}

\keycomment{\#<}{>\#}

\makeatother

\providecommand{\lemmaname}{Lemma}
\providecommand{\propositionname}{Proposition}
\providecommand{\remarkname}{Remark}
\providecommand{\theoremname}{Theorem}

\begin{document}

\title{Optimization Methods for Designing Sequences with Low Autocorrelation
Sidelobes}

\author{Junxiao Song, Prabhu Babu, and Daniel P. Palomar, \IEEEmembership{Fellow, IEEE}%
\thanks{Junxiao Song, Prabhu Babu, and Daniel P. Palomar are with the Hong
Kong University of Science and Technology (HKUST), Hong Kong. E-mail:
\{jsong, eeprabhubabu, palomar\}@ust.hk.%
}}
\maketitle
\begin{abstract}
Unimodular sequences with low autocorrelations are desired in many
applications, especially in the area of radar and code-division multiple
access (CDMA). In this paper, we propose a new algorithm to design
unimodular sequences with low integrated sidelobe level (ISL), which
is a widely used measure of the goodness of a sequence's correlation
property. The algorithm falls into the general framework of majorization-minimization
(MM) algorithms and thus shares the monotonic property of such algorithms.
In addition, the algorithm can be implemented via fast Fourier transform
(FFT) operations and thus is computationally efficient. Furthermore,
after some modifications the algorithm can be adapted to incorporate
spectral constraints, which makes the design more flexible. Numerical
experiments show that the proposed algorithms outperform existing
algorithms in terms of both the quality of designed sequences and
the computational complexity.\end{abstract}

\begin{IEEEkeywords}
Unimodular sequences, autocorrelation, integrated sidelobe level,
majorization-minimization.
\end{IEEEkeywords}

\section{Introduction}

Since the 1950s, digital communications engineers have sought to identify
sequences whose autocorrelation sidelobes are collectively as low
as possible according to some suitable measure of \textquotedblleft goodness\textquotedblright .
Applications of sequences with low autocorrelation sidelobes range
from synchronization to code-division multiple access (CDMA) systems
and also radar systems \cite{turyn1968sequences,golomb2005signal}.
Low autocorrelation can improve the detection performance of weak
targets \cite{levanon2004radar} in range compression radar and it
is also desired for synchronization purposes in CDMA systems. Additionally,
due to the limitations of sequence generation hardware components
(such as the maximum signal amplitude clip of analog-to-digital converters
and power amplifiers), it is usually more desirable to transmit unimodular
(i.e., constant modulus) sequences to maximize the transmitted power
available in the system \cite{he2012waveform}.

The focus of this paper is on the design of unimodular sequences with
low (aperiodic or periodic) autocorrelation sidelobes. Let $\{x_{n}\}_{n=1}^{N}$
denote the complex unimodular (without loss of generality, we will
assume the modulus to be 1) sequence to be designed. The aperiodic
($r_{k}$) and periodic ($\hat{r}_{k}$) autocorrelations of $\{x_{n}\}_{n=1}^{N}$
are defined as
\[
\begin{aligned}r_{k} & =\sum_{n=1}^{N-k}x_{n}x_{n+k}^{*}=r_{-k}^{*},\, k=0,\ldots,N-1,\\
\hat{r}_{k} & =\sum_{n=1}^{N}x_{n}x_{(n+k)({\rm mod}\, N)}^{*}=\hat{r}_{-k}^{*},\, k=0,\ldots,N-1.
\end{aligned}
\]
Then the integrated sidelobe level (ISL) of the aperiodic autocorrelations
is defined as 
\begin{equation}
{\rm ISL}=\sum_{k=1}^{N-1}|r_{k}|^{2},\label{eq:1}
\end{equation}
which will be used to measure the goodness of the aperiodic autocorrelation
property of a sequence. The ISL metric for the periodic autocorrelations
can be defined similarly. Note that the ISL metric is highly related
to another important goodness measure: the merit factor (MF), which
is defined as the ratio of the central lobe energy to the total energy
of all other lobes by Golay in 1972 \cite{MeritFactor1972}, i.e.,
\begin{equation}
{\rm MF}=\frac{\left|r_{0}\right|^{2}}{2\sum_{k=1}^{N-1}|r_{k}|^{2}}=\frac{N^{2}}{2{\rm ISL}}.\label{eq:merit_factor}
\end{equation}

Owing to the practical importance and widespread applications of sequences
with good autocorrelation properties, in particular with low ISL values
or large MF values, a lot of effort has been devoted to identifying
these sequences via either analytical construction methods or computational
approaches in the literature. The studies focused more on the design
of binary sequences at the early stage, and extended to polyphase
sequences later; see \cite{barker1953group,MeritFactor1972,Binary2golay1977,Binary4mertens1996,Binary1Kocabas2003,Binary5jedwab2005survey,Binary3wang2008}
for the case of binary sequences and  \cite{FrankSequence1963,Golomb1993polyphase,borwein2005polyphase,Polycode2009,De_Maio2009,stoica2009new,ITROX2012}
for the case of polyphase or unimodular sequences. Some sequences
with good autocorrelation properties that can be constructed in closed-form
have been proposed in the literature, such as the Frank sequence \cite{FrankSequence1963}
and the Golomb sequence \cite{Golomb1993polyphase}. Computational
approaches, such as exhaustive search \cite{Binary4mertens1996},
evolutionary algorithms \cite{Binary1Kocabas2003}, heuristic search
\cite{Binary3wang2008} and stochastic optimization methods \cite{borwein2005polyphase,Polycode2009},
have also been suggested. However, these methods are generally not
capable of designing long sequences, say $N\sim10^{3}$ or larger,
due to the increasing computational complexity. Recently, an algorithm
named CAN (cyclic algorithm new) was proposed in \cite{stoica2009new}
that can be used to produce unimodular sequences of length $N\sim10^{6}$
or even larger with low aperiodic ISL and it has been adapted to generate
sequences with impulse-like periodic autocorrelations \cite{PeCAN}
and incorporate spectral constraints \cite{he2010waveform}. But rather
than the original ISL metric, the CAN algorithm tries to minimize
another simpler criterion which is stated to be ``almost equivalent''
to the ISL metric.

In this paper, we develop an efficient algorithm that directly minimizes
the ISL metric (of both aperiodic and periodic autocorrelations) monotonically,
which we call MISL (monotonic minimizer for integrated sidelobe level).
The MISL algorithm is derived via applying the general majorization-minimization
(MM) method to the ISL minimization problem twice and admits a simple
closed-form solution at every iteration. The convergence of the algorithm
to a stationary point is proved. Similar to the CAN algorithm, the
proposed MISL algorithm can also be implemented via fast Fourier transform
(FFT) operations and is thus very efficient. But due to the nature
of double majorization, the MISL algorithm may converge very slow
especially for large $N,$ and we propose to apply some acceleration
schemes to fix this issue. A MISL-like algorithm, which shares the
simplicity and efficiency of MISL, has also been developed for minimizing
the ISL metric when there are additional spectral constraints.

The remaining sections of the paper are organized as follows. In Section
II, the problem formulations are presented and the CAN algorithm is
reviewed. In Section III, we first give a brief review of the MM framework
and then the MISL algorithm is derived, followed by the convergence
analysis. In Section IV, some acceleration schemes for MISL are presented
and then the spectral-MISL algorithm is developed in Section V. Finally,
Section VI presents some numerical examples and the conclusions are
given in Section VII.

\emph{Notation}: Boldface upper case letters denote matrices, boldface
lower case letters denote column vectors, and italics denote scalars.
$\mathbf{R}$ and $\mathbf{C}$ denote the real field and the complex
field, respectively. $\mathrm{Re}(\cdot)$ and $\mathrm{Im}(\cdot)$
denote the real and imaginary part respectively. ${\rm arg}(\cdot)$
denotes the phase of a complex number. The superscripts $(\cdot)^{T}$,
$(\cdot)^{*}$ and $(\cdot)^{H}$ denote transpose, complex conjugate,
and conjugate transpose, respectively. $X_{i,j}$ denotes the (\emph{i}-th,
\emph{j}-th) element of matrix $\mathbf{X}$ and $x_{i}$ denotes
the \emph{i}-th element of vector $\mathbf{x}$. $\mathrm{Tr}(\cdot)$
denotes the trace of a matrix. ${\rm diag}(\mathbf{X})$ is a column
vector consisting of all the diagonal elements of $\mathbf{X}$. ${\rm Diag}(\mathbf{x})$
is a diagonal matrix formed with $\mathbf{x}$ as its principal diagonal.
${\rm vec}(\mathbf{X})$ is a column vector consisting of all the
columns of $\mathbf{X}$ stacked. $\mathbf{I}_{n}$ denotes an $n\times n$
identity matrix.

\section{Problem Formulation and Existing Methods}

The problem of interest is the design of a complex unimodular sequence
$\{x_{n}\}_{n=1}^{N}$ that minimizes the ISL metric (of either aperiodic
or periodic autocorrelations), i.e., 
\begin{equation}
\begin{array}{ll}
\underset{x_{n}}{\mathsf{minimize}} & {\rm ISL}\\
\mathsf{subject\; to} & \left|x_{n}\right|=1,\, n=1,\ldots,N.
\end{array}\label{eq:ISL-minimize}
\end{equation}
In the following, we first reformulate the problem \eqref{eq:ISL-minimize}
for the case of aperiodic and periodic autocorrelations respectively
and then briefly review two corresponding algorithms proposed in the
literature.

\subsection{The Aperiodic Autocorrelation}

It has been shown in \cite{stoica2009new} that the ISL metric of
the aperiodic autocorrelations can be expressed in the frequency domain
as 
\begin{equation}
{\rm ISL}=\frac{1}{4N}\sum_{p=1}^{2N}\left[\left|\sum_{n=1}^{N}x_{n}e^{-j\omega_{p}(n-1)}\right|^{2}-N\right]^{2},\label{eq:2}
\end{equation}
where $\omega_{p}=\frac{2\pi}{2N}(p-1),\: p=1,\cdots,2N$. Thus, the
sequence design problem \eqref{eq:ISL-minimize} becomes 
\begin{equation}
\begin{array}{ll}
\underset{x_{n}}{\mathsf{minimize}} & {\displaystyle \sum_{p=1}^{2N}}\left[\left|{\displaystyle \sum_{n=1}^{N}}x_{n}e^{-j\omega_{p}(n-1)}\right|^{2}-N\right]^{2}\\
\mathsf{subject\; to} & \left|x_{n}\right|=1,\, n=1,\ldots,N.
\end{array}\label{eq:ISL-freq}
\end{equation}

Let us define 
\begin{eqnarray}
\mathbf{x} & = & [x_{1},\cdots,x_{N}]^{T}\\
\mathbf{a}_{p} & = & [1,e^{j\omega_{p}},\cdots,e^{j\omega_{p}(N-1)}]^{T},\, p=1,\ldots,2N,
\end{eqnarray}
then the problem \eqref{eq:ISL-freq} can be rewritten as 
\begin{equation}
\begin{array}{ll}
\underset{\mathbf{x}}{\mathsf{minimize}} & {\displaystyle \sum_{p=1}^{2N}}\left[\mathbf{a}_{p}^{H}\mathbf{x}\mathbf{x}^{H}\mathbf{a}_{p}-N\right]^{2}\\
\mathsf{subject\; to} & \left|x_{n}\right|=1,\, n=1,\ldots,N.
\end{array}\label{eq:prob-fre}
\end{equation}
Expanding the square in the objective yields 
\begin{equation}
\begin{array}{ll}
\underset{\mathbf{x}}{\mathsf{minimize}} & {\displaystyle \sum_{p=1}^{2N}}\left(\left(\mathbf{a}_{p}^{H}\mathbf{x}\mathbf{x}^{H}\mathbf{a}_{p}\right)^{2}-2N\mathbf{a}_{p}^{H}\mathbf{x}\mathbf{x}^{H}\mathbf{a}_{p}+N^{2}\right)\\
\mathsf{subject\; to} & \left|x_{n}\right|=1,\, n=1,\ldots,N,
\end{array}
\end{equation}
where the second term in the objective can be shown to be a constant
(using Parseval's theorem), i.e., $\sum_{p=1}^{2N}\left|\mathbf{a}_{p}^{H}\mathbf{x}\right|^{2}=2N\left\Vert \mathbf{x}\right\Vert _{2}^{2}=2N^{2}.$
Thus by ignoring the constant terms, the problem can be simplified
as 
\begin{equation}
\begin{array}{ll}
\underset{\mathbf{x}}{\mathsf{minimize}} & {\displaystyle \sum_{p=1}^{2N}}\left(\mathbf{a}_{p}^{H}\mathbf{x}\mathbf{x}^{H}\mathbf{a}_{p}\right)^{2}\\
\mathsf{subject\; to} & \left|x_{n}\right|=1,\, n=1,\ldots,N.
\end{array}\label{eq:prob-simple}
\end{equation}
The problem is hard to tackle, due to the quartic objective function
and also the nonconvex unit-modulus constraints.

\subsection{The Periodic Autocorrelation}

Similarly, in the periodic case, the ISL metric ${\rm ISL}=\sum_{k=1}^{N-1}|\hat{r}_{k}|^{2}$
can also be expressed in the frequency domain as \cite{PeCAN}
\begin{equation}
{\rm ISL}=\frac{1}{N}\sum_{p=1}^{N}\left[\left|\sum_{n=1}^{N}x_{n}e^{-j\hat{\omega}_{p}(n-1)}\right|^{2}-N\right]^{2},\label{eq:ISL-pe}
\end{equation}
where $\hat{\omega}_{p}=\frac{2\pi}{N}(p-1),\: p=1,\cdots,N$. By
defining 
\begin{eqnarray}
\hat{\mathbf{a}}_{p} & = & [1,e^{j\hat{\omega}_{p}},\cdots,e^{j\hat{\omega}_{p}(N-1)}]^{T},\, p=1,\ldots,N,\label{eq:a_p_hat}
\end{eqnarray}
the ISL minimization problem \eqref{eq:ISL-minimize} in this case
can be rewritten as 
\begin{equation}
\begin{array}{ll}
\underset{\mathbf{x}}{\mathsf{minimize}} & {\displaystyle \sum_{p=1}^{N}}\left[\hat{\mathbf{a}}_{p}^{H}\mathbf{x}\mathbf{x}^{H}\hat{\mathbf{a}}_{p}-N\right]^{2}\\
\mathsf{subject\; to} & \left|x_{n}\right|=1,\, n=1,\ldots,N.
\end{array}\label{eq:prob-fre-pe}
\end{equation}
Following similar steps as in the previous subsection, the problem
can be further simplified as 
\begin{equation}
\begin{array}{ll}
\underset{\mathbf{x}}{\mathsf{minimize}} & {\displaystyle \sum_{p=1}^{N}}\left(\hat{\mathbf{a}}_{p}^{H}\mathbf{x}\mathbf{x}^{H}\hat{\mathbf{a}}_{p}\right)^{2}\\
\mathsf{subject\; to} & \left|x_{n}\right|=1,\, n=1,\ldots,N.
\end{array}\label{eq:prob-simple-pe}
\end{equation}

We note that although problem \eqref{eq:prob-simple-pe} looks similar
to the problem \eqref{eq:prob-simple} in the aperiodic case, it is
usually a simpler task to find sequences with good periodic autocorrelation
properties. In particular, in the periodic case unimodular sequences
with zero ISL exist for any length $N$ \cite{PeCAN}, however it
is not possible to construct unimodular sequences with exact impulse-like
aperiodic autocorrelations.

\subsection{Existing Methods}

An algorithm named CAN was proposed in \cite{stoica2009new} for designing
unimodular sequences with low ISL of aperiodic autocorrelations. But
instead of minimizing the ISL metric directly, i.e., solving \eqref{eq:ISL-freq},
the authors proposed to solve the following simpler problem
\begin{equation}
\begin{array}{ll}
\underset{x_{n},\psi_{p}}{\mathsf{minimize}} & {\displaystyle \sum_{p=1}^{2N}}\left|{\displaystyle \sum_{n=1}^{N}}x_{n}e^{-j\omega_{p}(n-1)}-\sqrt{N}e^{j\psi_{p}}\right|^{2}\\
\mathsf{subject\; to} & \left|x_{n}\right|=1,\, n=1,\ldots,N,
\end{array}\label{eq:prob-CAN}
\end{equation}
whose objective is quadratic in $\{x_{n}\}$. Let $\mathbf{A}$ be
the following $N\times2N$ matrix: 
\begin{equation}
\mathbf{A}=[\mathbf{a}_{1},\ldots,\mathbf{a}_{2N}],\label{eq:Amat}
\end{equation}
then the objective function of problem \eqref{eq:prob-CAN} can be
rewritten in the following more compact form: 
\begin{equation}
\left\Vert \mathbf{A}^{H}\mathbf{x}-\sqrt{N}\mathbf{v}\right\Vert ^{2},\label{eq:obj_CAN_compact}
\end{equation}
where $\mathbf{v}=[e^{j\psi_{1}},\ldots,e^{j\psi_{2N}}]^{T}.$ The
CAN algorithm then minimizes \eqref{eq:obj_CAN_compact} by alternating
between $\mathbf{x}$ and $\mathbf{v}$. For a given $\mathbf{x}$,
the minimization of \eqref{eq:obj_CAN_compact} with respect to $\{\psi_{p}\}$
yields 
\begin{equation}
\psi_{p}={\rm arg}(f_{p}),\, p=1,\ldots,2N,
\end{equation}
where $\mathbf{f}=\mathbf{A}^{H}\mathbf{x}$. Similarly, for a given
$\mathbf{v}$, the minimizing $\mathbf{x}$ is given by 
\begin{equation}
x_{n}=e^{j{\rm arg}(g_{n})},\, n=1,\ldots,N,
\end{equation}
where $\mathbf{g}=\mathbf{A}\mathbf{v}.$ The CAN algorithm is summarized
in Algorithm \ref{alg:CAN} for ease of comparison later.

\begin{algorithm}[tbh]
\begin{algor}[1]
\item [{Require:}] sequence length $N$
\item [{{*}}] Set $k=0$, initialize $\mathbf{x}^{(0)}$. 
\item [{repeat}]~

\begin{algor}[1]
\item [{{*}}] $\mathbf{f}=\mathbf{A}^{H}\mathbf{x}^{(k)}$
\item [{{*}}] $v_{p}=e^{j{\rm arg}(f_{p})},\, p=1,\ldots,2N$
\item [{{*}}] $\mathbf{g}=\mathbf{A}\mathbf{v}$
\item [{{*}}] $x_{n}^{(k+1)}=e^{j{\rm arg}(g_{n})},\, n=1,\ldots,N$
\item [{{*}}] $k\leftarrow k+1$
\end{algor}
\item [{until}] convergence
\end{algor}
\protect\caption{\label{alg:CAN} The CAN algorithm proposed in \cite{stoica2009new}.}
\end{algorithm}

It was mentioned in \cite{stoica2009new} that the problem \eqref{eq:prob-CAN}
is ``almost equivalent'' to the original problem \eqref{eq:ISL-freq}.
But as the authors also pointed out, the two problems are not exactly
equivalent and they may have different solutions. It is easy to see
that minimizing with respect to $\{\psi_{p}\}$ in \eqref{eq:prob-CAN}
leads to
\begin{equation}
\begin{array}{ll}
\underset{x_{n}}{\mathsf{minimize}} & {\displaystyle \sum_{p=1}^{2N}}\left[\left|{\displaystyle \sum_{n=1}^{N}}x_{n}e^{-j\omega_{p}(n-1)}\right|-\sqrt{N}\right]^{2}\\
\mathsf{subject\; to} & \left|x_{n}\right|=1,\, n=1,\ldots,N,
\end{array}
\end{equation}
which is different from the original problem \eqref{eq:ISL-freq}.
So the point the CAN algorithm converges to is not necessary a local
minimum (or even a stationary point) of the original ISL metric minimization
problem. Although CAN does not directly minimize the ISL metric, in
\cite{stoica2009new} the authors showed numerically that CAN could
generate sequences with good correlation properties. Moreover, CAN
can be implemented via FFT, thus can be used to design very long sequences.

In the periodic case, an algorithm called PeCAN was used in \cite{PeCAN}
to construct unimodular sequences with low ISL. It can be viewed as
an adaptation of the CAN algorithm for the periodic case and similarly
instead of minimizing the original ISL metric \eqref{eq:ISL-pe},
the PeCAN algorithm considered the simpler criterion 
\begin{equation}
\sum_{p=1}^{N}\left|\sum_{n=1}^{N}x_{n}e^{-j\hat{\omega}_{p}(n-1)}-\sqrt{N}e^{j\psi_{p}}\right|^{2},
\end{equation}
where $\{\psi_{p}\}$ are auxiliary variables as in \eqref{eq:prob-CAN}.
It was shown numerically in \cite{PeCAN} that PeCAN can generate
almost perfect (with zero ISL) sequences from random initializations.

Although the criteria considered by CAN and PeCAN are somehow related
to the original ISL metric and the performance of the two algorithms
is acceptable in some situations, one may expect that directly solving
the original ISL minimization problem can probably lead to better
performance. In the next section, we will develop an algorithm that
directly minimizes the ISL metric, while at the same time is computationally
as efficient as the CAN algorithm.

\section{ISL Minimization via Majorization-Minimization}

In this section, we first introduce the general majorization-minimization
(MM) method briefly and then develop a simple algorithm for the ISL
minimization problem based on the general MM scheme. In the derivation,
we will focus on the aperiodic case, i.e., the problem \eqref{eq:prob-simple},
but the resulting algorithm can be readily adapted for the periodic
case.

\subsection{The MM Method\label{sub:MM-Method}}

The MM method refers to the majorization-minimization method, which
is an approach to solve optimization problems that are too difficult
to solve directly. The principle behind the MM method is to transform
a difficult problem into a series of simple problems. Interested readers
may refer to \cite{hunter2004MMtutorial} and references therein for
more details (recent generalizations include \cite{razaviyayn2013unified,Aldo2013}).

Suppose we want to minimize $f(\mathbf{x})$ over $\mathcal{X}\in\mathbf{C}^{n}$.
Instead of minimizing the cost function $f(\mathbf{x})$ directly,
the MM approach optimizes a sequence of approximate objective functions
that majorize $f(\mathbf{x})$. More specifically, starting from a
feasible point $\mathbf{x}^{(0)},$ the algorithm produces a sequence
$\{\mathbf{x}^{(k)}\}$ according to the following update rule 
\begin{equation}
\mathbf{x}^{(k+1)}\in\underset{\mathbf{x}\in\mathcal{X}}{\arg\min}\,\, u(\mathbf{x},\mathbf{x}^{(k)}),\label{eq:major_update}
\end{equation}
where $\mathbf{x}^{(k)}$ is the point generated by the algorithm
at iteration $k,$ and $u(\mathbf{x},\mathbf{x}^{(k)})$ is the majorization
function of $f(\mathbf{x})$ at $\mathbf{x}^{(k)}$. Formally, the
function $u(\mathbf{x},\mathbf{x}^{(k)})$ is said to majorize the
function $f(\mathbf{x})$ at the point $\mathbf{x}^{(k)}$ provided
\begin{eqnarray}
u(\mathbf{x},\mathbf{x}^{(k)}) & \geq & f(\mathbf{x}),\quad\forall\mathbf{x}\in\mathcal{X},\label{eq:major1}\\
u(\mathbf{x}^{(k)},\mathbf{x}^{(k)}) & = & f(\mathbf{x}^{(k)}).\label{eq:major2}
\end{eqnarray}
In other words, function $u(\mathbf{x},\mathbf{x}^{(k)})$ is an upper
bound of $f(\mathbf{x})$ over $\mathcal{X}$ and coincides with $f(\mathbf{x})$
at $\mathbf{x}^{(k)}$.

To summarize, to minimize $f(\mathbf{x})$ over $\mathcal{X}\in\mathbf{C}^{n}$,
the main steps of the majorization-minimization scheme are 
\begin{enumerate}
\item Find a feasible point $\mathbf{x}^{(0)}$ and set $k=0.$ 
\item Construct a majorization function $u(\mathbf{x},\mathbf{x}^{(k)})$
of $f(\mathbf{x})$ at $\mathbf{x}^{(k)}$.
\item Let $\mathbf{x}^{(k+1)}\in\underset{\mathbf{x}\in\mathcal{X}}{\arg\min}\,\, u(\mathbf{x},\mathbf{x}^{(k)}).$ 
\item If some convergence criterion is met, exit; otherwise, set $k=k+1$
and go to step (2). 
\end{enumerate}
It is easy to show that with this scheme, the objective value is decreased
monotonically at every iteration, i.e., 
\begin{equation}
f(\mathbf{x}^{(k+1)})\leq u(\mathbf{x}^{(k+1)},\mathbf{x}^{(k)})\leq u(\mathbf{x}^{(k)},\mathbf{x}^{(k)})=f(\mathbf{x}^{(k)}).\label{eq:descent-property}
\end{equation}
The first inequality and the third equality follow from the the properties
of the majorization function, namely \eqref{eq:major1} and \eqref{eq:major2}
respectively and the second inequality follows from \eqref{eq:major_update}.
The monotonicity makes MM algorithms very stable in practice.

\subsection{MISL}

To solve the aperiodic problem \eqref{eq:prob-simple} via majorization-minimization,
the key step is to find a majorization function of the objective such
that the majorized problem is easy to solve. For that purpose we first
present a simple result that will be useful when constructing the
majorization function. 
\begin{lem}
\label{lem:majorizer}Let $\mathbf{L}$ be an $n\times n$ Hermitian
matrix and $\mathbf{M}$ be another $n\times n$ Hermitian matrix
such that $\mathbf{M}\succeq\mathbf{L}.$\textup{ }Then for any point
$\mathbf{x}_{0}\in\mathbf{C}^{n}$, the quadratic function $\mathbf{x}^{H}\mathbf{L}\mathbf{x}$
is majorized by $\mathbf{x}^{H}\mathbf{M}\mathbf{x}+2{\rm Re}\left(\mathbf{x}^{H}(\mathbf{L}-\mathbf{M})\mathbf{x}_{0}\right)+\mathbf{x}_{0}^{H}(\mathbf{M}-\mathbf{L})\mathbf{x}_{0}$
at\textup{ $\mathbf{x}_{0}$.}\end{lem}
\begin{IEEEproof}
It is easy to check that the two functions are equal at point $\mathbf{x}_{0}$.
Since $\mathbf{M}\succeq\mathbf{L},$ we further have 
\[
\begin{aligned} & \,\,\mathbf{x}^{H}\mathbf{L}\mathbf{x}\\
= & \,\,\mathbf{x}_{0}^{H}\mathbf{L}\mathbf{x}_{0}\!+\!2{\rm Re}\left((\mathbf{x}-\mathbf{x}_{0})^{H}\mathbf{L}\mathbf{x}_{0}\right)\!+\!(\mathbf{x}-\mathbf{x}_{0})^{H}\mathbf{L}(\mathbf{x}-\mathbf{x}_{0})\\
\leq & \,\,\mathbf{x}_{0}^{H}\mathbf{L}\mathbf{x}_{0}\!+\!2{\rm Re}\left((\mathbf{x}-\mathbf{x}_{0})^{H}\mathbf{L}\mathbf{x}_{0}\right)\!+\!(\mathbf{x}-\mathbf{x}_{0})^{H}\mathbf{M}(\mathbf{x}-\mathbf{x}_{0})\\
= & \,\,\mathbf{x}^{H}\mathbf{M}\mathbf{x}\!+\!2{\rm Re}\left(\mathbf{x}^{H}(\mathbf{L}-\mathbf{M})\mathbf{x}_{0}\right)\!+\!\mathbf{x}_{0}^{H}(\mathbf{M}-\mathbf{L})\mathbf{x}_{0}
\end{aligned}
\]
for any $\mathbf{x}\in\mathbf{C}^{n}$. The proof is complete. 
\end{IEEEproof}
The objective of the problem \eqref{eq:prob-simple} is quartic with
respect to $\mathbf{x}$. To find a majorization function, some reformulations
are necessary. Let us define $\mathbf{X}=\mathbf{x}\mathbf{x}^{H}$
and $\mathbf{A}_{p}=\mathbf{a}_{p}\mathbf{a}_{p}^{H}$, then the problem
\eqref{eq:prob-simple} can be rewritten as 
\begin{equation}
\begin{array}{ll}
\underset{\mathbf{x},\mathbf{X}}{\mathsf{minimize}} & {\displaystyle \sum_{p=1}^{2N}}{\rm Tr}(\mathbf{X}\mathbf{A}_{p})^{2}\\
\mathsf{subject\; to} & \mathbf{X}=\mathbf{x}\mathbf{x}^{H}\\
 & \left|x_{n}\right|=1,\, n=1,\ldots,N.
\end{array}\label{eq:prob-capitalX}
\end{equation}
Since ${\rm Tr}(\mathbf{X}\mathbf{A}_{p})={\rm vec}(\mathbf{X})^{H}{\rm vec}(\mathbf{A}_{p})$,
problem \eqref{eq:prob-capitalX} is just 
\begin{equation}
\begin{array}{ll}
\underset{\mathbf{x},\mathbf{X}}{\mathsf{minimize}} & {\rm vec}(\mathbf{X})^{H}\mathbf{\Phi}{\rm vec}(\mathbf{X})\\
\mathsf{subject\; to} & \mathbf{X}=\mathbf{x}\mathbf{x}^{H}\\
 & \left|x_{n}\right|=1,\, n=1,\ldots,N,
\end{array}\label{eq:prob-vecX}
\end{equation}
where $\mathbf{\Phi}=\sum_{p=1}^{2N}{\rm vec}(\mathbf{A}_{p}){\rm vec}(\mathbf{A}_{p})^{H}.$
We can see that the objective of \eqref{eq:prob-vecX} is a quadratic
function in $\mathbf{X}$ now. Given $\mathbf{X}^{(k)}=\mathbf{x}^{(k)}\left(\mathbf{x}^{(k)}\right)^{H}$
at iteration $k$, by choosing $\mathbf{M}=\lambda_{{\rm max}}(\mathbf{\Phi})\mathbf{I}$
in Lemma \ref{lem:majorizer} we know that the objective of \eqref{eq:prob-vecX}
is majorized by the following function at $\mathbf{X}^{(k)}$:
\begin{equation}
\begin{aligned} & u_{1}(\mathbf{X},\mathbf{X}^{(k)})\\
= & \lambda_{{\rm max}}(\mathbf{\Phi}){\rm vec}(\mathbf{X})^{H}{\rm vec}(\mathbf{X})\\
 & +2{\rm Re}\left({\rm vec}(\mathbf{X})^{H}(\mathbf{\Phi}-\lambda_{{\rm max}}(\mathbf{\Phi})\mathbf{I}){\rm vec}(\mathbf{X}^{(k)})\right)\\
 & +{\rm vec}(\mathbf{X}^{(k)})^{H}(\lambda_{{\rm max}}(\mathbf{\Phi})\mathbf{I}-\mathbf{\Phi}){\rm vec}(\mathbf{X}^{(k)})
\end{aligned}
\label{eq:major_vecX}
\end{equation}
where $\lambda_{{\rm max}}(\mathbf{\Phi})$ is the maximum eigenvalue
of $\mathbf{\Phi}$ and can be shown to be $2N^{2}$ (see Appendix
\ref{sub:Proof-Lambda_max}). Since ${\rm vec}(\mathbf{X})^{H}{\rm vec}(\mathbf{X})=(\mathbf{x}^{H}\mathbf{x})^{2}=N^{2}$,
the first term of \eqref{eq:major_vecX} is just a constant. After
ignoring the constant terms in \eqref{eq:major_vecX}, the majorized
problem of \eqref{eq:prob-vecX} is given by
\begin{equation}
\begin{array}{ll}
\underset{\mathbf{x},\mathbf{X}}{\mathsf{minimize}} & {\rm Re}\left({\rm vec}(\mathbf{X})^{H}\left(\mathbf{\Phi}-2N^{2}\mathbf{I}\right){\rm vec}(\mathbf{X}^{(k)})\right)\\
\mathsf{subject\; to} & \mathbf{X}=\mathbf{x}\mathbf{x}^{H}\\
 & \left|x_{n}\right|=1,\, n=1,\ldots,N,
\end{array}
\end{equation}
which can be rewritten as
\begin{equation}
\begin{array}{ll}
\underset{\mathbf{x},\mathbf{X}}{\mathsf{minimize}} & {\displaystyle \sum_{p=1}^{2N}}{\rm Tr}(\mathbf{X}^{(k)}\mathbf{A}_{p}){\rm Tr}(\mathbf{A}_{p}\mathbf{X})-2N^{2}{\rm Tr}(\mathbf{X}^{(k)}\mathbf{X})\\
\mathsf{subject\; to} & \mathbf{X}=\mathbf{x}\mathbf{x}^{H}\\
 & \left|x_{n}\right|=1,\, n=1,\ldots,N.
\end{array}\label{eq:MM-iterk}
\end{equation}

\begin{rem}
To solve the problem \eqref{eq:MM-iterk}, one possible approach is
to apply an SDP relaxation, yielding 
\begin{equation}
\begin{array}{ll}
\underset{\mathbf{X}}{\mathsf{minimize}} & {\displaystyle \sum_{p=1}^{2N}}{\rm Tr}(\mathbf{X}^{(k)}\mathbf{A}_{p}){\rm Tr}(\mathbf{A}_{p}\mathbf{X})-2N^{2}{\rm Tr}(\mathbf{X}^{(k)}\mathbf{X})\\
\mathsf{subject\; to} & \mathbf{X}\succeq\mathbf{0}\\
 & {\rm diag}(\mathbf{X})=\mathbf{1}.
\end{array}
\end{equation}
The above SDP is empirically observed to always give a rank one solution,
but proving the tightness of the SDR is out of the scope of this paper.
In addition, the computational complexity of solving an SDP is high
and it needs to be solved at every iteration, which makes it not amenable
for large $N.$ 
\end{rem}
In the following, we propose to apply majorization-minimization to
the problem \eqref{eq:MM-iterk} again, which will yield a simple
algorithm. To apply the second majorization, we first rewrite the
problem \eqref{eq:MM-iterk} as 
\begin{equation}
\begin{array}{ll}
\underset{\mathbf{x}}{\mathsf{minimize}} & {\displaystyle \sum_{p=1}^{2N}}\left|\mathbf{a}_{p}^{H}\mathbf{x}^{(k)}\right|^{2}\left|\mathbf{a}_{p}^{H}\mathbf{x}\right|^{2}-2N^{2}\left|\mathbf{x}^{H}\mathbf{x}^{(k)}\right|^{2}\\
\mathsf{subject\; to} & \left|x_{n}\right|=1,\, n=1,\ldots,N.
\end{array}
\end{equation}
It can be written more compactly as 
\begin{equation}
\begin{array}{ll}
\underset{\mathbf{x}}{\mathsf{minimize}} & \mathbf{x}^{H}\left(\mathbf{A}{\rm Diag}(\mathbf{p}^{(k)})\mathbf{A}^{H}-2N^{2}\mathbf{x}^{(k)}\left(\mathbf{x}^{(k)}\right)^{H}\right)\mathbf{x}\\
\mathsf{subject\; to} & \left|x_{n}\right|=1,\, n=1,\ldots,N,
\end{array}\label{eq:MM-2}
\end{equation}
where $\mathbf{A}=[\mathbf{a}_{1},\ldots,\mathbf{a}_{2N}]$ and $\mathbf{p}^{(k)}=\left|\mathbf{A}^{H}\mathbf{x}^{(k)}\right|^{2}$
($\left|\cdot\right|^{2}$ denotes the element-wise absolute-squared
value). We can clearly see that the objective function in \eqref{eq:MM-2}
is quadratic in $\mathbf{x}$ and by choosing $\mathbf{M}=p_{{\rm max}}^{(k)}\mathbf{A}\mathbf{A}^{H}$
in Lemma \ref{lem:majorizer} it is majorized by the following function
at $\mathbf{x}^{(k)}$: 
\begin{equation}
\begin{aligned} & u_{2}(\mathbf{x},\mathbf{x}^{(k)})\\
= & p_{{\rm max}}^{(k)}\mathbf{x}^{H}\mathbf{A}\mathbf{A}^{H}\mathbf{x}+2{\rm Re}\left(\mathbf{x}^{H}(\tilde{\mathbf{A}}-2N^{2}\mathbf{x}^{(k)}(\mathbf{x}^{(k)})^{H})\mathbf{x}^{(k)}\right)\\
 & +(\mathbf{x}^{(k)})^{H}(2N^{2}\mathbf{x}^{(k)}(\mathbf{x}^{(k)})^{H}-\tilde{\mathbf{A}})\mathbf{x}^{(k)}
\end{aligned}
\label{eq:major_x}
\end{equation}
where $p_{{\rm max}}^{(k)}=\max_{p}\{p_{p}^{(k)}:p=1,\ldots,2N\}$
and $\tilde{\mathbf{A}}=\mathbf{A}\left({\rm Diag}(\mathbf{p}^{(k)})-p_{{\rm max}}^{(k)}\mathbf{I}\right)\mathbf{A}^{H}$
. Since $\mathbf{A}\mathbf{A}^{H}=2N\mathbf{I},$ the first term of
\eqref{eq:major_x} is a constant. By ignoring the constant terms
in \eqref{eq:major_x}, we have the majorized problem of \eqref{eq:MM-2}:
\begin{equation}
\begin{array}{ll}
\underset{\mathbf{x}}{\mathsf{minimize}} & {\rm Re}\left(\mathbf{x}^{H}\left(\tilde{\mathbf{A}}-2N^{2}\mathbf{x}^{(k)}(\mathbf{x}^{(k)})^{H}\right)\mathbf{x}^{(k)}\right)\\
\mathsf{subject\; to} & \left|x_{n}\right|=1,\, n=1,\ldots,N,
\end{array}\label{eq:MM-prob2}
\end{equation}
which can be rewritten as 
\begin{equation}
\begin{array}{ll}
\underset{\mathbf{x}}{\mathsf{minimize}} & \left\Vert \mathbf{x}-\mathbf{y}\right\Vert _{2}\\
\mathsf{subject\; to} & \left|x_{n}\right|=1,\, n=1,\ldots,N,
\end{array}\label{eq:x-y}
\end{equation}
where
\begin{equation}
\begin{aligned}\mathbf{y} & =-\left(\tilde{\mathbf{A}}-2N^{2}\mathbf{x}^{(k)}\big(\mathbf{x}^{(k)}\big)^{H}\right)\mathbf{x}^{(k)}\\
 & =-\mathbf{A}\left({\rm Diag}(\mathbf{p}^{(k)})-p_{{\rm max}}^{(k)}\mathbf{I}-N^{2}\mathbf{I}\right)\mathbf{A}^{H}\mathbf{x}^{(k)}.
\end{aligned}
\label{eq:y_k}
\end{equation}
It is easy to see that the problem \eqref{eq:x-y} has a closed-form
solution, which is given by 
\begin{equation}
x_{n}=e^{j{\rm arg}(y_{n})},\, n=1,\ldots,N.
\end{equation}

Now we can summarize the overall algorithm and it is given in Algorithm
\ref{alg:MISL-1}. Note that although in the derivation we have applied
the majorization-minimization scheme twice, it can be viewed as directly
majorizing the objective of the problem \eqref{eq:prob-simple} at
$\mathbf{x}^{(k)}$ over the constraint set by the following function:
\begin{equation}
\begin{aligned} & u(\mathbf{x},\mathbf{x}^{(k)})\\
= & 2u_{2}(\mathbf{x},\mathbf{x}^{(k)})+4N^{4}-{\displaystyle \sum_{p=1}^{2N}}\left|\mathbf{a}_{p}^{H}\mathbf{x}^{(k)}\right|^{4}\\
= & 4{\rm Re}\left(\mathbf{x}^{H}\mathbf{A}\left({\rm Diag}(\mathbf{p}^{(k)})-p_{{\rm max}}^{(k)}\mathbf{I}-N^{2}\mathbf{I}\right)\mathbf{A}^{H}\mathbf{x}^{(k)}\right)\\
 & +8N^{2}\left(p_{{\rm max}}^{(k)}+N^{2}\right)-3{\displaystyle \sum_{p=1}^{2N}}\left|\mathbf{a}_{p}^{H}\mathbf{x}^{(k)}\right|^{4}
\end{aligned}
\label{eq:upper_u}
\end{equation}
and the algorithm preserves the monotonicity of the general majorization-minimization
method. Thus, we call the algorithm MISL (Monotonic minimizer for
Integrated Sidelobe Level). 

We also note that the above derivations can be carried out similarly
for the periodic problem \eqref{eq:prob-simple-pe} and it turns out
that the MISL algorithm can be easily applied for the periodic case
after replacing $\mathbf{A}$ in \eqref{eq:Amat} and Algorithm \ref{alg:MISL-1}
with $\hat{\mathbf{A}}=[\hat{\mathbf{a}}_{1},\ldots,\hat{\mathbf{a}}_{N}]$,
where $\{\hat{\mathbf{a}}_{p}\}_{p=1}^{N}$ are defined in \eqref{eq:a_p_hat}.

\begin{algorithm}[tbh]
\begin{algor}[1]
\item [{Require:}] sequence length $N$
\item [{{*}}] Set $k=0$, initialize $\mathbf{x}^{(0)}$. 
\item [{repeat}]~

\begin{algor}[1]
\item [{{*}}] $\mathbf{p}^{(k)}=\left|\mathbf{A}^{H}\mathbf{x}^{(k)}\right|^{2}$
\item [{{*}}] $p_{{\rm max}}^{(k)}=\max_{p}\{p_{p}^{(k)}:p=1,\ldots,2N\}$
\item [{{*}}] \begin{raggedright}
$\mathbf{y}=-\mathbf{A}\left({\rm Diag}(\mathbf{p}^{(k)})-p_{{\rm max}}^{(k)}\mathbf{I}-N^{2}\mathbf{I}\right)\mathbf{A}^{H}\mathbf{x}^{(k)}$
\par\end{raggedright}
\item [{{*}}] $x_{n}^{(k+1)}=e^{j{\rm arg}(y_{n})},\, n=1,\ldots,N$
\item [{{*}}] $k\leftarrow k+1$
\end{algor}
\item [{until}] convergence
\end{algor}
\protect\caption{\label{alg:MISL-1}MISL - Monotonic minimizer for Integrated Sidelobe
Level.}
\end{algorithm}

\subsection{Convergence Analysis}

The MISL algorithm given in Algorithm \ref{alg:MISL-1} is based on
the general maximization-minimization scheme, thus according to subsection
\ref{sub:MM-Method}, we know that the sequence of objective values
evaluated at $\{\mathbf{x}^{(k)}\}$ generated by the algorithm is
nonincreasing. And it is easy to see that the objective value is bounded
below by $0$, thus the sequence of objective values is guaranteed
to converge to a finite value.

In this section, we will further analyze the convergence property
of the sequence $\{\mathbf{x}^{(k)}\}$ generated by the MISL algorithm
and show the convergence to a stationary point.

We first introduce a first-order optimality condition for minimizing
a smooth function over an arbitrary constraint set, which follows
from \cite{Bertsekas2003}.
\begin{prop}
Let $f:\mathbf{R}^{n}\rightarrow\mathbf{R}$ be a smooth function,
and let $\mathbf{x}^{\star}$ be a local minimum of $f$ over a subset
$\mathcal{X}$ of $\mathbf{R}^{n}$. Then 
\begin{equation}
\nabla f(\mathbf{x}^{\star})^{T}\mathbf{y}\geq0,\,\forall\mathbf{y}\in T_{\mathcal{X}}(\mathbf{x}^{\star}),\label{eq:opt_condition}
\end{equation}
where $T_{\mathcal{X}}(\mathbf{x}^{\star})$ denotes the tangent cone
of $\mathcal{X}$ at $\mathbf{x}^{\star}.$
\end{prop}
A point $\mathbf{x}\in\mathcal{X}$ satisfying the first-order optimality
condition \eqref{eq:opt_condition} will be referred as a stationary
point.

Next we exploit a nice property shared by the minimization problems
under consideration.
\begin{lem}
\label{lem:prob-equivalent}Let $g:\mathbf{C}^{N}\rightarrow\mathbf{R}$
be a real valued function with complex variables. Then the problem
\begin{equation}
\begin{array}{ll}
\underset{\mathbf{x}\in\mathbf{C}^{N}}{\mathsf{minimize}} & g(\mathbf{x})\\
\mathsf{subject\; to} & \left|x_{n}\right|=1,\, n=1,\ldots,N
\end{array}\label{eq:prob_eq1}
\end{equation}
is equivalent to 
\begin{equation}
\begin{array}{ll}
\underset{\mathbf{x}\in\mathbf{C}^{N}}{\mathsf{minimize}} & g(\mathbf{x})+\alpha\mathbf{x}^{H}\mathbf{x}\\
\mathsf{subject\; to} & \left|x_{n}\right|=1,\, n=1,\ldots,N
\end{array}\label{eq:prob_eq2}
\end{equation}
for any finite value $\alpha$, in the sense that they share the same
optimal solutions.\end{lem}
\begin{IEEEproof}
Since the additional term $\alpha\mathbf{x}^{H}\mathbf{x}$ is just
the constant $\alpha N$ over the constraint set $\{\mathbf{x}\in\mathbf{C}^{N}|\left|x_{n}\right|=1,n=1,\ldots,N\}$,
it is obvious that any solution of problem \eqref{eq:prob_eq1} is
also a solution of problem \eqref{eq:prob_eq2} and vice versa.
\end{IEEEproof}
To facilitate the analysis, we further note that upon defining 

\begin{equation}
\tilde{\mathbf{x}}=[{\rm Re}(\mathbf{x})^{T},{\rm Im}(\mathbf{x})^{T}]^{T}\label{eq:x_real}
\end{equation}
\begin{equation}
\tilde{\mathbf{A}}_{p}=\left[\begin{array}{cc}
{\rm Re}(\mathbf{A}_{p}) & -{\rm Im}(\mathbf{A}_{p})\\
{\rm Im}(\mathbf{A}_{p}) & {\rm Re}(\mathbf{A}_{p})
\end{array}\right],\label{eq:Ap_real}
\end{equation}
it is straightforward to show that the quartic complex minimization
problem \eqref{eq:prob-simple} is equivalent to the following real
one:
\begin{equation}
\begin{array}{ll}
\underset{\tilde{\mathbf{x}}}{\mathsf{minimize}} & {\displaystyle \sum_{p=1}^{2N}}\left(\tilde{\mathbf{x}}^{T}\tilde{\mathbf{A}}_{p}\tilde{\mathbf{x}}\right)^{2}\\
\mathsf{subject\; to} & \tilde{x}_{n}^{2}+\tilde{x}_{n+N}^{2}=1,\, n=1,\ldots,N.
\end{array}\label{eq:prob_real_equivalent}
\end{equation}

We are now ready to state the convergence properties of MISL.
\begin{thm}
Let $\{\mathbf{x}^{(k)}\}$ be the sequence generated by the MISL
algorithm in Algorithm \ref{alg:MISL-1}. Then every limit point of
the sequence $\{\mathbf{x}^{(k)}\}$ is a stationary point of the
problem \eqref{eq:prob-simple}.\end{thm}
\begin{IEEEproof}
Denote the objective functions of the problem \eqref{eq:prob-simple}
and its real equivalent \eqref{eq:prob_real_equivalent} by $f(\mathbf{x})$
and $\tilde{f}(\tilde{\mathbf{x}})$, respectively. Denote the constraint
sets of the problem \eqref{eq:prob-simple} and \eqref{eq:prob_real_equivalent}
by $\mathcal{C}$ and $\tilde{\mathcal{C}}$, respectively, i.e.,
$\mathcal{C}=\{\mathbf{x}\in\mathbf{C}^{N}|\left|x_{n}\right|=1,n=1,\ldots,N\}$
and $\tilde{\mathcal{C}}=\{\tilde{\mathbf{x}}\in\mathbf{R}^{2N}|\tilde{x}_{n}^{2}+\tilde{x}_{n+N}^{2}=1,\, n=1,\ldots,N\}$.
Then from the derivation of MISL, we know that $f(\mathbf{x})$ is
majorized by the function $u(\mathbf{x},\mathbf{x}^{(k)})$ in \eqref{eq:upper_u}
at $\mathbf{x}^{(k)}$ over $\mathcal{C}$. According to the general
MM scheme described in subsection \ref{sub:MM-Method}, we have
\[
f(\mathbf{x}^{(k+1)})\leq u(\mathbf{x}^{(k+1)},\mathbf{x}^{(k)})\leq u(\mathbf{x}^{(k)},\mathbf{x}^{(k)})=f(\mathbf{x}^{(k)}),
\]
which means $\{f(\mathbf{x}^{(k)})\}$ is a nonincreasing sequence.

Assume that there exists a converging subsequence $\mathbf{x}^{(k_{j})}\rightarrow\mathbf{x}^{(\infty)},$
then 
\[
\begin{aligned}u(\mathbf{x}^{(k_{j+1})},\mathbf{x}^{(k_{j+1})}) & =f(\mathbf{x}^{(k_{j+1})})\leq f(\mathbf{x}^{(k_{j}+1)})\\
 & \leq u(\mathbf{x}^{(k_{j}+1)},\mathbf{x}^{(k_{j})})\leq u(\mathbf{x},\mathbf{x}^{(k_{j})}),\forall\mathbf{x}\in\mathcal{C}.
\end{aligned}
\]
Letting $j\rightarrow+\infty,$ we obtain 
\begin{equation}
u(\mathbf{x}^{(\infty)},\mathbf{x}^{(\infty)})\leq u(\mathbf{x},\mathbf{x}^{(\infty)}),\,\forall\mathbf{x}\in\mathcal{C},
\end{equation}
i.e., $\mathbf{x}^{(\infty)}$ is a global minimizer of $u(\mathbf{x},\mathbf{x}^{(\infty)})$
over $\mathcal{C}$. According to Lemma \ref{lem:prob-equivalent},
$\mathbf{x}^{(\infty)}$ is also a global minimizer of
\begin{equation}
\begin{array}{ll}
\underset{\mathbf{x}}{\mathsf{minimize}} & u(\mathbf{x},\mathbf{x}^{(\infty)})+4N^{3}\mathbf{x}^{H}\mathbf{x}\\
\mathsf{subject\; to} & \mathbf{x}\in\mathcal{C}.
\end{array}\label{eq:prob-mm}
\end{equation}
With the definitions of $\tilde{\mathbf{x}}$ and $\tilde{\mathbf{A}}_{p}$
given in \eqref{eq:x_real} and \eqref{eq:Ap_real} and by ignoring
some constant terms in $u(\mathbf{x},\mathbf{x}^{(\infty)})$, it
is easy to show that \eqref{eq:prob-mm} is equivalent to 
\begin{equation}
\begin{array}{ll}
\underset{\tilde{\mathbf{x}}}{\mathsf{minimize}} & 2p_{{\rm max}}^{(\infty)}{\displaystyle \sum_{p=1}^{2N}}\tilde{\mathbf{x}}^{T}\tilde{\mathbf{A}}_{p}\tilde{\mathbf{x}}+4\tilde{\mathbf{x}}^{T}\mathbf{b}+4N^{3}\tilde{\mathbf{x}}^{T}\tilde{\mathbf{x}}\\
\mathsf{subject\; to} & \tilde{\mathbf{x}}\in\tilde{\mathcal{C}},
\end{array}\label{eq:prob-mm-real}
\end{equation}
where 
\begin{equation}
\mathbf{b}=\sum_{p=1}^{2N}\left((\tilde{\mathbf{x}}^{(\infty)})^{T}\tilde{\mathbf{A}}_{p}\tilde{\mathbf{x}}^{(\infty)}-p_{{\rm max}}^{(\infty)}\right)\tilde{\mathbf{A}}_{p}\tilde{\mathbf{x}}^{(\infty)}-2N^{3}\tilde{\mathbf{x}}^{(\infty)}.
\end{equation}
Thus, $\tilde{\mathbf{x}}^{(\infty)}=[{\rm Re}(\mathbf{x}^{(\infty)})^{T},{\rm Im}(\mathbf{x}^{(\infty)})^{T}]^{T}$
is a global minimizer of \eqref{eq:prob-mm-real}. Then as a necessary
condition, we have 
\begin{equation}
\nabla\tilde{u}(\tilde{\mathbf{x}}^{(\infty)})^{T}\mathbf{y}\geq0,\,\forall\mathbf{y}\in T_{\tilde{\mathcal{C}}}(\tilde{\mathbf{x}}^{(\infty)}),
\end{equation}
where $\tilde{u}(\tilde{\mathbf{x}})$ denotes the objective function
of \eqref{eq:prob-mm-real}. It is easy to check that 
\[
\nabla\tilde{f}(\tilde{\mathbf{x}}^{(\infty)})=\nabla\tilde{u}(\tilde{\mathbf{x}}^{(\infty)})=4\sum_{p=1}^{2N}\left((\tilde{\mathbf{x}}^{(\infty)})^{T}\tilde{\mathbf{A}}_{p}\tilde{\mathbf{x}}^{(\infty)}\right)\tilde{\mathbf{A}}_{p}\tilde{\mathbf{x}}^{(\infty)}.
\]
Thus we have 
\begin{equation}
\nabla\tilde{f}(\tilde{\mathbf{x}}^{(\infty)})^{T}\mathbf{y}\geq0,\,\forall\mathbf{y}\in T_{\tilde{\mathcal{C}}}(\tilde{\mathbf{x}}^{(\infty)}),
\end{equation}
implying that $\tilde{\mathbf{x}}^{(\infty)}$ is a stationary point
of the problem \eqref{eq:prob_real_equivalent}. Due to the equivalence
of problem \eqref{eq:prob_real_equivalent} and \eqref{eq:prob-simple},
the proof is complete.
\end{IEEEproof}

\subsection{Computational Complexity of MISL}

From Algorithm \ref{alg:MISL-1}, we can see that the per iteration
computational complexity of MISL is dominated by two matrix vector
multiplications involving $\mathbf{A}$, which can be easily computed
via FFT and IFFT operations. More specifically, let $\mathbf{F}$
be the $2N\times2N$ DFT matrix with $F_{m,n}=e^{-j\frac{2mn\pi}{2N}},0\leq m,n<2N$,
it is easy to see that the $N\times2N$ matrix $\mathbf{A}$ is composed
of the first $N$ rows of $\mathbf{F}^{H}$. Given a vector $\mathbf{x}\in\mathbf{C}^{N}$,
the product $\mathbf{A}^{H}\mathbf{x}$ can be computed via the FFT
of the zero padded version of $\mathbf{x}$, i.e., $\mathbf{F}[\mathbf{x}^{H},\mathbf{0}_{N}^{H}]^{H}$.
Similarly, given a vector $\mathbf{z}\in\mathbf{C}^{2N}$, the product
\textbf{$\mathbf{A}\mathbf{z}$} consists of the first $N$ elements
of the inverse FFT of $\mathbf{z}$ (i.e., $\mathbf{F}^{H}\mathbf{z}$).
In the periodic case, since the matrix $\hat{\mathbf{A}}^{H}$ itself
is the $N\times N$ DFT matrix, multiplying a vector by $\hat{\mathbf{A}}^{H}$
and $\hat{\mathbf{A}}$ is just the FFT and IFFT of the vector. Thus,
the MISL algorithm is computationally very efficient and can be used
for the design of very long sequences, say $N\sim10^{6}$.

\section{Acceleration Schemes \label{sec:Acceleration-Schemes}}

As described in the previous section, the derivation of MISL is based
on majorization-minimization principle, and the nature of the majorization
functions will dictate the convergence speed of the algorithm. From
numerical simulations, we noticed that for large $N$, the convergence
of MISL is very slow, which may be due to the double majorization
scheme that we carried out in its derivation. One option to fix this
issue is to employ some acceleration schemes to accelerate the convergence
of MISL. In the following, we will consider two such schemes (for
the aperiodic case) and the resulting algorithms can be adapted for
the periodic case by using $\hat{\mathbf{A}}$ instead of $\mathbf{A}$.

\subsection{Acceleration via Fixed Point Theory \label{sub:Acceleration-SQUAREM}}

In this subsection, we describe an acceleration scheme that can be
applied to accelerate the MISL algorithm proposed in this paper. It
is the so called squared iterative method (SQUAREM), which was originally
proposed in \cite{SQUAREM} to accelerate any EM algorithms. SQUAREM
follows the idea of the Cauchy-Barzilai-Borwein (CBB) method \cite{raydan2002CBB},
which combines the classical steepest descent method and the two-point
step size gradient method \cite{barzilai1988_BB}, to solve the nonlinear
fixed-point problem of EM. It only requires the EM updating scheme
and can be readily implemented as an 'off-the-shelf' accelerator.
Since MM is a generalization of EM and the update rule is just a fixed-point
iteration, SQUAREM can be readily applied to MM algorithms.

Let $\mathbf{F}_{{\rm MISL}}(\cdot)$ denote the nonlinear fixed-point
iteration map of the MISL algorithm: 
\begin{equation}
\mathbf{x}^{(k+1)}=\mathbf{F_{{\rm MISL}}}(\mathbf{x}^{(k)}),\label{eq:acc1}
\end{equation}
whose form can be defined by the following equation: 
\begin{equation}
\mathbf{x}^{(k+1)}=e^{j{\rm arg}\left(-\mathbf{A}\left({\rm Diag}(\mathbf{p}^{(k)})-p_{{\rm max}}^{(k)}\mathbf{I}-N^{2}\mathbf{I}\right)\mathbf{A}^{H}\mathbf{x}^{(k)}\right)},\label{eq:acc_update}
\end{equation}
where $\mathbf{p}^{(k)}$ and $p_{{\rm max}}^{(k)}$ are the same
as in Algorithm \ref{alg:MISL-1} and functions $e^{(\cdot)}$ and
${\rm arg(\cdot)}$ are applied element-wise to the vectors. With
this, the steps of the accelerated-MISL based on SQUAREM are summarized
in Algorithm \ref{alg:acc_MISL-1}.

A problem of the general SQUAREM is that it may violate the nonlinear
constraints, so in Algorithm \ref{alg:acc_MISL-1} the function $e^{j{\rm arg}(\cdot)}$
has been applied to project wayward points back to the feasible region.
A second problem of SQUAREM is that it can violate the descent property
of the original MM algorithm. To ensure the descent property, a strategy
based on backtracking is adopted, which repeatedly halves the distance
between $\alpha$ and $-1$ (i.e., $\alpha\leftarrow(\alpha-1)/2$)
until the descent property is maintained. To see why this works, we
first note that $\mathbf{x}=e^{j{\rm arg}(\mathbf{x}^{(k)}-2\alpha\mathbf{r}+\alpha^{2}\mathbf{v})}=\mathbf{x}_{2}$
if $\alpha=-1.$ In addition, since ${\rm ISL}(\mathbf{x}_{2})\leq{\rm ISL}(\mathbf{x}^{(k)})$
due to the descent property of original MM steps, ${\rm ISL}(\mathbf{x})\leq{\rm ISL}(\mathbf{x}^{(k)})$
is guaranteed to hold as $\alpha\rightarrow-1.$ In practice, usually
only a few backtracking steps are needed to maintain the monotonicity
of the algorithm.

\begin{algorithm}[tbh]
\begin{algor}[1]
\item [{Require:}] sequence length $N$ 
\item [{{*}}] Set $k=0$, initialize $\mathbf{x}^{(0)}$. 
\item [{repeat}]~

\begin{algor}[1]
\item [{{*}}] $\mathbf{x}_{1}=\mathbf{F_{{\rm MISL}}}(\mathbf{x}^{(k)})$ 
\item [{{*}}] $\mathbf{x}_{2}=\mathbf{F_{{\rm MISL}}}(\mathbf{x}_{1})$ 
\item [{{*}}] $\mathbf{r}=\mathbf{x}_{1}-\mathbf{x}^{(k)}$ 
\item [{{*}}] $\mathbf{v}=\mathbf{x}_{2}-\mathbf{x}_{1}-\mathbf{r}$ 
\item [{{*}}] Compute the step-length $\alpha=-\frac{\|\mathbf{r}\|}{\|\mathbf{v}\|}$ 
\item [{{*}}] $\mathbf{x}=e^{j{\rm arg}(\mathbf{x}^{(k)}-2\alpha\mathbf{r}+\alpha^{2}\mathbf{v})}$ 
\item [{while}] ${\rm ISL}(\mathbf{x})>{\rm ISL}(\mathbf{x}^{(k)})$
\item [{{*}}] $\alpha\leftarrow(\alpha-1)/2$
\item [{{*}}] $\mathbf{x}=e^{j{\rm arg}(\mathbf{x}^{(k)}-2\alpha\mathbf{r}+\alpha^{2}\mathbf{v})}$ 
\item [{endwhile}]~
\item [{{*}}] $\mathbf{x}^{(k+1)}=\mathbf{x}$
\item [{{*}}] $k\leftarrow k+1$ 
\end{algor}
\item [{until}] convergence 
\end{algor}
\protect\caption{\label{alg:acc_MISL-1}Accelerated-MISL.}
\end{algorithm}

\subsection{Acceleration via Backtracking}

As mentioned earlier, the slow convergence of MISL could be due to
the double majorization, which may have resulted in a majorization
function that is a very loose approximation of the original objective
function. Thus to accelerate the convergence, one possibility is to
find a better approximation of the objective at each iteration. Notice
that to preserve the descent property \eqref{eq:descent-property}
of the majorization-minimization method, it only requires $u(\mathbf{x},\mathbf{x}^{(k)})\geq f(\mathbf{x})$
at $\mathbf{x}=\mathbf{x}^{(k+1)},$ i.e., at the minimizer of $u(\mathbf{x},\mathbf{x}^{(k)})$,
rather than requiring $u(\mathbf{x},\mathbf{x}^{(k)})$ to be a global
upper bound of $f(\mathbf{x}).$ With this in mind, we consider the
following approximation function of $f(\mathbf{x})=\sum_{p=1}^{2N}\left|\mathbf{a}_{p}^{H}\mathbf{x}\right|^{4}$
at $\mathbf{x}^{(k)}$: 

\begin{equation}
\begin{aligned}u_{L}(\mathbf{x},\mathbf{x}^{(k)})= & \,4{\rm Re}\left(\mathbf{x}^{H}\mathbf{A}\left({\rm Diag}(\mathbf{p}^{(k)})-L\mathbf{I}\right)\mathbf{A}^{H}\mathbf{x}^{(k)}\right)\\
 & +8N^{2}L-3{\displaystyle \sum_{p=1}^{2N}}\left|\mathbf{a}_{p}^{H}\mathbf{x}^{(k)}\right|^{4},
\end{aligned}
\end{equation}
where $L\geq0$ needs to be chosen such that $u_{L}(\mathbf{x},\mathbf{x}^{(k)})\geq f(\mathbf{x})$
at the minimizer of $u_{L}(\mathbf{x},\mathbf{x}^{(k)})$ over the
constraint set $\mathcal{C}=\{\mathbf{x}\in\mathbf{C}^{N}|\left|x_{n}\right|=1,n=1,\ldots,N\}.$
It is easy to see that the minimizer of $u_{L}(\mathbf{x},\mathbf{x}^{(k)})$
over $\mathcal{C}$ is given by 
\begin{equation}
\mathbf{x}_{L}^{\star}=e^{j{\rm arg}\left(\mathbf{A}\left(L\mathbf{I}-{\rm Diag}(\mathbf{p}^{(k)})\right)\mathbf{A}^{H}\mathbf{x}^{(k)}\right)},
\end{equation}
where functions $e^{(\cdot)}$ and ${\rm arg(\cdot)}$ are applied
element-wise to the vectors. Thus we need to choose $L$ such that
\begin{equation}
u_{L}(\mathbf{x}_{L}^{\star},\mathbf{x}^{(k)})\geq f(\mathbf{x}_{L}^{\star}).\label{eq:inequal_at_minimum}
\end{equation}
It is worth noting that by taking $L=p_{{\rm max}}^{(k)}+N^{2}$,
the function $u_{L}(\mathbf{x},\mathbf{x}^{(k)})$ is just the majorization
function of $f(\mathbf{x})$ in \eqref{eq:upper_u}. Hence we know
that the condition \eqref{eq:inequal_at_minimum} can be guaranteed
if $L$ is greater than $p_{{\rm max}}^{(k)}+N^{2}$. To find a better
approximation, here we adopt a backtracking strategy to choose $L$.
More specifically, at iteration $k$ we choose $L$ to be the smallest
element in $\{p_{{\rm max}}^{(k)}+(2^{i_{k}}-1)N\}_{i_{k}=0,1,\ldots}$
such that condition \eqref{eq:inequal_at_minimum} is satisfied. Note
that by choosing $L$ in this way, the resulting function $u_{L}(\mathbf{x},\mathbf{x}^{(k)})$
is not guaranteed to be a global upper bound anymore, but the monotonic
property of the algorithm is ensured. We call this algorithm as backtracking-MISL
and the steps of the algorithm are summarized in Algorithm \ref{alg:MISL-heuristic}.

\begin{algorithm}[tbh]
\begin{algor}[1]
\item [{Require:}] sequence length $N$ 
\item [{{*}}] Set $k=0$, initialize $\mathbf{x}^{(0)}$. 
\item [{repeat}]~

\begin{algor}[1]
\item [{{*}}] $\mathbf{p}^{(k)}=\left|\mathbf{A}^{H}\mathbf{x}^{(k)}\right|^{2}$
\item [{{*}}] $p_{{\rm max}}^{(k)}=\max_{p}\{p_{p}^{(k)}:p=1,\ldots,2N\}$ 
\item [{{*}}] Set $i_{k}=0$
\item [{repeat}]~
\item [{{*}}] $L=p_{{\rm max}}^{(k)}+(2^{i_{k}}-1)N$
\item [{{*}}] $\mathbf{x}_{L}^{\star}=e^{j{\rm arg}\left(\mathbf{A}\left(L\mathbf{I}-{\rm Diag}(\mathbf{p}^{(k)})\right)\mathbf{A}^{H}\mathbf{x}^{(k)}\right)}$
\item [{{*}}] $i_{k}\leftarrow i_{k}+1$
\item [{until}] $u_{L}(\mathbf{x}_{L}^{\star},\mathbf{x}^{(k)})\geq f(\mathbf{x}_{L}^{\star})$
\item [{{*}}] $\mathbf{x}^{(k+1)}=\mathbf{x}_{L}^{\star}$
\item [{{*}}] $k\leftarrow k+1$ 
\end{algor}
\item [{until}] convergence 
\end{algor}
\protect\caption{\label{alg:MISL-heuristic}Backtracking-MISL.}
\end{algorithm}

\section{ISL Minimization with spectral constraints }

In some applications, for example in cognitive radar \cite{cognitiveRadar_Haykin2006},
apart from designing sequences with good correlation properties they
need to satisfy some spectral constraints like the spectral power
in some frequency bands should be lower than some specified level.
In practice, designing such sequences is a challenge, there are few
methods available in the literature dealing with this problem, see
\cite{he2010waveform}. In this section, we will show how MISL can
be adapted to design sequences with spectral constraints and we call
the resulting algorithm spectral-MISL.

The spectral constraints such as the power in some band, denoted hereafter
by the set of indices $\Omega\subset[1,2N]$, should be lower than
a threshold can be expressed as:
\begin{equation}
\sum_{k\in\Omega}|\mathbf{a}_{k}^{H}\mathbf{x}|^{2}\:\leq\epsilon,\label{eq:SMISL1}
\end{equation}
where $\epsilon$ denotes some pre-specified threshold. Then the problem
of minimizing ISL subject to some spectral constraints can be expressed
as
\begin{equation}
\begin{array}{ll}
\underset{\mathbf{x}}{\mathsf{minimize}} & {\displaystyle \sum_{k=1}^{N-1}}|r_{k}|^{2}\\
\mathsf{subject\; to} & \left|x_{n}\right|=1,\, n=1,\ldots,N\\
 & \sum_{k\in\Omega}^ {}|\mathbf{a}_{k}^{H}\mathbf{x}|^{2}\:\leq\epsilon.
\end{array}\label{eq:SMISL2}
\end{equation}
For a given $\epsilon>0$, we can always find a $\lambda$ such that
the problem \eqref{eq:SMISL2} can be transformed into the following
equivalent problem:
\begin{equation}
\begin{array}{ll}
\underset{\mathbf{x}}{\mathsf{minimize}} & {\displaystyle \sum_{k=1}^{N-1}}|r_{k}|^{2}+\lambda{\displaystyle \sum_{k\in\Omega}}|\mathbf{a}_{k}^{H}\mathbf{x}|^{2}\\
\mathsf{subject\; to} & \left|x_{n}\right|=1,\, n=1,\ldots,N.
\end{array}\label{eq:SMISL3}
\end{equation}
As shown before in \eqref{eq:prob-fre}, the problem can be expressed
as
\begin{equation}
\begin{array}{ll}
\underset{\mathbf{x}}{\mathsf{minimize}} & {\displaystyle \sum_{p=1}^{2N}\left[\mathbf{a}_{p}^{H}\mathbf{x}\mathbf{x}^{H}\mathbf{a}_{p}-N\right]^{2}+\lambda\sum_{k\in\Omega}^ {}|\mathbf{a}_{k}^{H}\mathbf{x}|^{2}}\\
\mathsf{subject\; to} & \left|x_{n}\right|=1,\, n=1,\ldots,N.
\end{array}\label{eq:SMISL4}
\end{equation}
From here on we can follow the derivation of MISL to derive the spectral-MISL
algorithm for problem \eqref{eq:SMISL4}. For instance, after the
first majorization at a given point $\mathbf{x}^{(k)}$, the majorized
problem is 
\begin{equation}
\begin{array}{ll}
\underset{\mathbf{x}}{\mathsf{minimize}} & 2{\displaystyle \sum_{p=1}^{2N}}\left|\mathbf{a}_{p}^{H}\mathbf{x}^{(k)}\right|^{2}\mathbf{x}^{H}\mathbf{a}_{p}\mathbf{a}_{p}^{H}\mathbf{x}+\lambda{\displaystyle \sum_{k\in\Omega}}|\mathbf{a}_{k}^{H}\mathbf{x}|^{2}\\
 & -4N^{2}\mathbf{x}^{H}\mathbf{x}^{(k)}(\mathbf{x}^{(k)})^{H}\mathbf{x}\\
\mathsf{subject\; to} & \left|x_{n}\right|=1,\, n=1,\ldots,N,
\end{array}\label{eq:SMISL5}
\end{equation}
which can be rewritten as
\begin{equation}
\begin{array}{ll}
\underset{\mathbf{x}}{\mathsf{minimize}} & 2{\displaystyle \sum_{p=1}^{2N}}\bar{p}_{p}^{(k)}\mathbf{x}^{H}\mathbf{a}_{p}\mathbf{a}_{p}^{H}\mathbf{x}-4N^{2}\mathbf{x}^{H}\mathbf{x}^{(k)}(\mathbf{x}^{(k)})^{H}\mathbf{x}\\
\mathsf{subject\; to} & \left|x_{n}\right|=1,\, n=1,\ldots,N,
\end{array}\label{eq:SMISL6}
\end{equation}
where
\begin{equation}
\bar{p}_{p}^{(k)}=\begin{cases}
\left|\mathbf{a}_{p}^{H}\mathbf{x}^{(k)}\right|^{2}+\lambda/2, & p\in\Omega\\
\left|\mathbf{a}_{p}^{H}\mathbf{x}^{(k)}\right|^{2}, & {\rm otherwise}.
\end{cases}\label{eq:SMISL6_2}
\end{equation}
By defining $\bar{\mathbf{p}}^{(k)}=\left[\bar{p}_{1}^{(k)},\ldots,\bar{p}_{2N}^{(k)}\right]^{T}$
and $\mathbf{A}=[\mathbf{a}_{1},\ldots,\mathbf{a}_{2N}]$, we get
the problem
\begin{equation}
\begin{array}{ll}
\underset{\mathbf{x}}{\mathsf{minimize}} & \mathbf{x}^{H}\left(\mathbf{A}{\rm Diag}(\bar{\mathbf{p}}^{(k)})\mathbf{A}^{H}-2N^{2}\mathbf{x}^{(k)}\left(\mathbf{x}^{(k)}\right)^{H}\right)\mathbf{x}\\
\mathsf{subject\; to} & \left|x_{n}\right|=1,\, n=1,\ldots,N,
\end{array}\label{eq:SMISL7}
\end{equation}
which has the same form as the problem \eqref{eq:MM-2}. So we can
follow the same steps as before and derive the spectral-MISL algorithm,
the main steps of which are listed in Algorithm \ref{alg:SMISL}.
The acceleration schemes described in section \ref{sec:Acceleration-Schemes}
can be readily applied to spectral-MISL as well, however we will not
elaborate on that here.

\begin{algorithm}[tbh]
\begin{algor}[1]
\item [{Require:}] sequence length $N$, index set $\Omega$ and $\lambda$. 
\item [{{*}}] Set $k=0$, initialize $\mathbf{x}^{(0)}$. 
\item [{repeat}]~

\begin{algor}[1]
\item [{{*}}] \begin{raggedright}
$\bar{p}_{p}^{(k)}=\begin{cases}
\left|\mathbf{a}_{p}^{H}\mathbf{x}^{(k)}\right|^{2}+\lambda/2, & p\in\Omega\\
\left|\mathbf{a}_{p}^{H}\mathbf{x}^{(k)}\right|^{2}, & {\rm otherwise}
\end{cases}$ 
\par\end{raggedright}
\item [{{*}}] $\bar{p}_{{\rm max}}^{(k)}=\max_{p}\{\bar{p}_{p}^{(k)}:p=1,\ldots,2N\}$ 
\item [{{*}}] \begin{raggedright}
$\mathbf{y}=-\mathbf{A}\left({\rm Diag}(\mathbf{\bar{p}}^{(k)})-\bar{p}_{{\rm max}}^{(k)}\mathbf{I}-N^{2}\mathbf{I}\right)\mathbf{A}^{H}\mathbf{x}^{(k)}$ 
\par\end{raggedright}
\item [{{*}}] $x_{n}^{(k+1)}=e^{j{\rm arg}(y_{n})},\, n=1,\ldots,N$ 
\item [{{*}}] $k\leftarrow k+1$ 
\end{algor}
\item [{until}] convergence 
\end{algor}
\protect\caption{\label{alg:SMISL} Spectral-MISL.}
\end{algorithm}

\section{Numerical Experiments}

To compare the performance of the proposed MISL algorithms with existing
ones and to show the potential of MISL to design sequences for various
scenarios, we present some experimental results in this section. All
experiments were performed on a PC with a 3.20GHz i5-3470 CPU and
8GB RAM.

\subsection{ISL Minimization}

In this subsection, we present numerical results on applying the proposed
algorithms to design unimodular sequences with low ISL and compare
the performance with the CAN algorithm \cite{stoica2009new}, which
is known to be computationally very efficient. The Matlab code of
CAN was downloaded from the website%
\footnote{http://www.sal.ufl.edu/book/%
} of the book \cite{he2012waveform}.

\begin{figure}[t]
\centering{}\includegraphics[width=0.95\columnwidth]{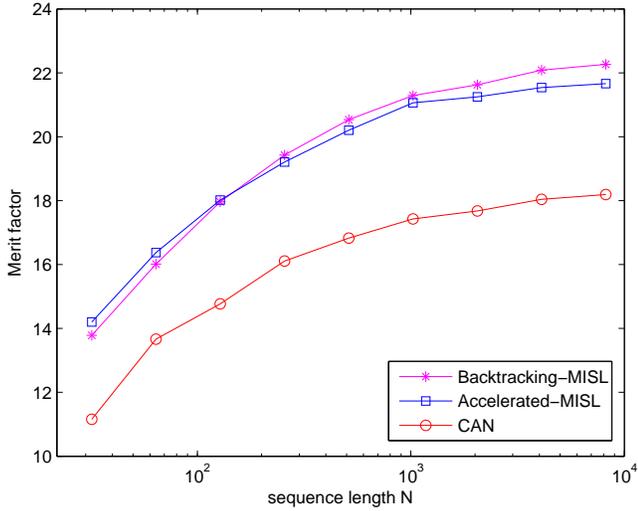}
\protect\caption{\label{fig:avg_MF}Merit factor (MF) versus sequence length. Each
curve is averaged over 100 random trials.}
\end{figure}

\begin{figure}[t]
\centering{}\includegraphics[width=0.95\columnwidth]{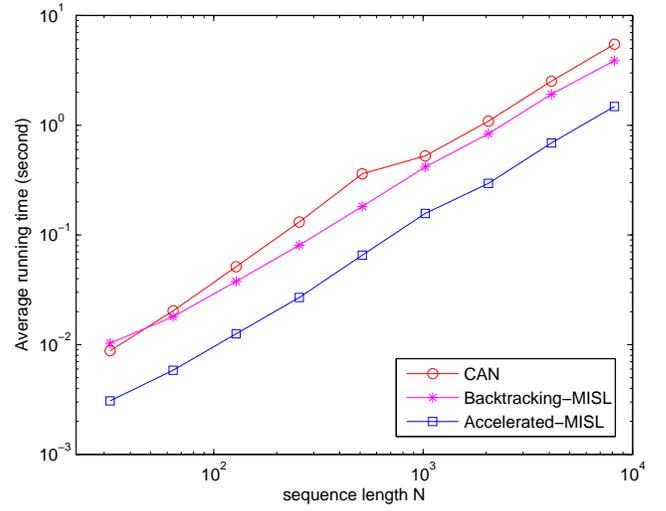}
\protect\caption{\label{fig:avg_run_time}Average running time versus sequence length.
Each curve is averaged over 100 random trials.}
\end{figure}

\begin{figure}[t]
\centering{}\includegraphics[width=0.95\columnwidth]{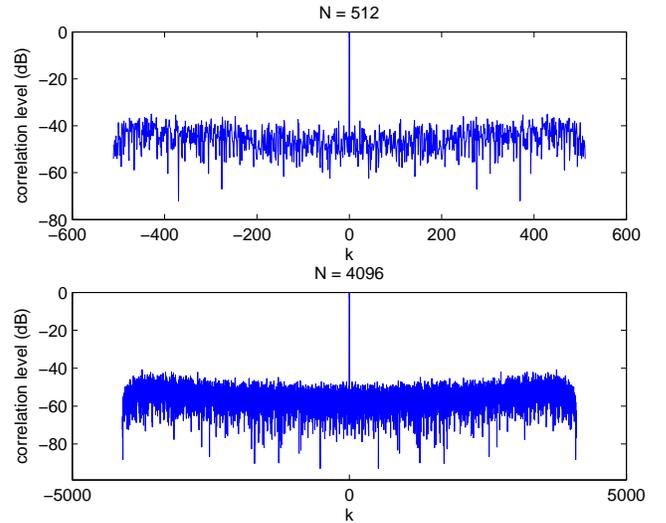}
\protect\caption{\label{fig:MISL_correlation_level}Correlation level of two sequences
of length $N=512$ and $N=4096$, generated by the accelerated-MISL
algorithm.}
\end{figure}

We first compare the quality, measured by the merit factor (MF) defined
in \eqref{eq:merit_factor} (recall that the larger the better), of
the output sequences generated by different algorithms. In this experiment,
for all algorithms, the stopping criterion was set to be $\left|{\rm ISL}(\mathbf{x}^{(k+1)})-{\rm ISL}(\mathbf{x}^{(k)})\right|/\max\left(1,{\rm ISL}(\mathbf{x}^{(k)})\right)\leq10^{-5}$
and the initial sequence $\{x_{n}^{(0)}\}_{n=1}^{N}$ was chosen to
be $\{e^{j2\pi\theta_{n}}\}_{n=1}^{N}$, where $\{\theta_{n}\}_{n=1}^{N}$
are independent random variables uniformly distributed in $[0,1]$.
Each algorithm was repeated 100 times for each of the following lengths:
$N=2^{5},2^{6},\ldots,2^{13}$. The average merit factor of the output
sequences and the average running time are shown in Figs. \ref{fig:avg_MF}
and \ref{fig:avg_run_time} for different algorithms. From Fig. \ref{fig:avg_MF},
we can see that the proposed backtracking-MISL and accelerated-MISL
algorithms can generate sequences with consistently larger MF than
the CAN algorithm for all lengths; and at the same time they are faster
than the CAN algorithm as can be seen from Fig. \ref{fig:avg_run_time},
especially the accelerated-MISL algorithm. Taking both the merit factor
and the computational complexity into consideration, the accelerated-MISL
seems to be the best among the three algorithms. The correlation levels
of two example sequences of length $N=512$ and $N=4096$ generated
by the accelerated-MISL algorithm are shown in Fig. \ref{fig:MISL_correlation_level},
where the correlation level is defined as 
\[
\textrm{correlation level}=20\log_{10}\left|\frac{r_{k}}{r_{0}}\right|,k=1-N,\ldots,N-1.
\]

\begin{figure}[htbp]
\begin{centering}
\includegraphics[width=0.95\columnwidth]{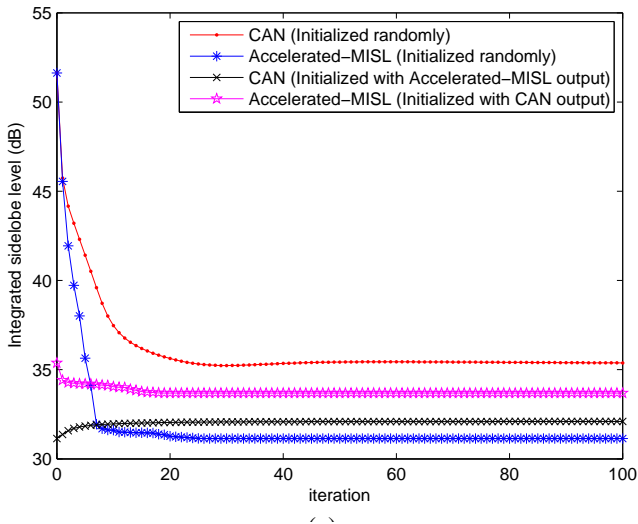}
\par\end{centering}

\begin{centering}
\begin{minipage}[t]{1\columnwidth}%
\begin{center}
(a)
\par\end{center}%
\end{minipage}
\par\end{centering}

\medskip{}

\begin{centering}
\includegraphics[width=0.95\columnwidth]{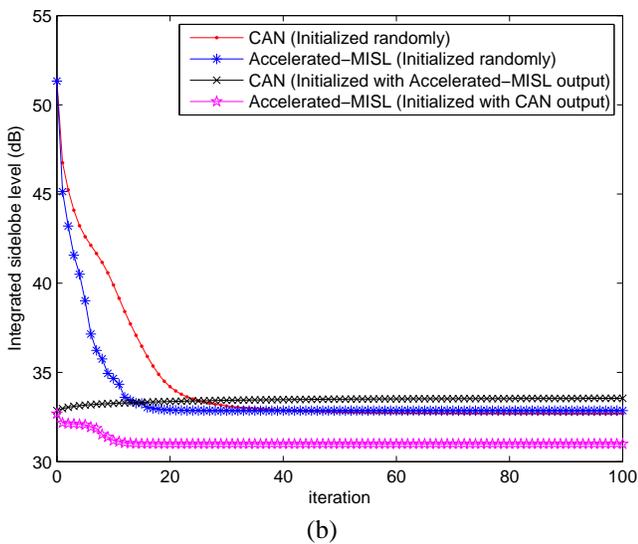}
\par\end{centering}

\begin{centering}
\begin{minipage}[t]{1\columnwidth}%
\begin{center}
(b)
\par\end{center}%
\end{minipage}
\par\end{centering}

\centering{}\protect\caption{\label{fig:convergence}Evolution of the Integrated sidelobe level
(ISL) in two different random trails with $N=32$.}
\end{figure}

Next, we perform an experiment to show the different convergence properties
of the proposed accelerated-MISL algorithm and the CAN algorithm.
For both algorithms, we first initialize them with the same randomly
generated sequence. Then we initialize CAN with the output sequence
of accelerated-MISL and initialize accelerated-MISL with the output
sequence of CAN, and run the two algorithms again. The evolution curves
of the ISL in two different random trails with $N=32$ are shown in
Fig. \ref{fig:convergence}(a) and \ref{fig:convergence}(b). From
the figures, we can see that when initialized with the output sequence
of CAN, the accelerated-MISL can further decrease the ISL a little
bit, which means the point CAN converges to is probably not a local
minimum. While in contrast, when initialized with the output sequence
of accelerated-MISL, CAN increases the ISL, which means the point
accelerated-MISL converges to is probably a local optimal and also
tells the fact that CAN does not share the monotonicity of MISL. In
addition, we can also see that when initialized with different sequences,
CAN and accelerated-MISL may converge to different points. Thus, initializing
the algorithms with good sequences, e.g. Frank sequence \cite{FrankSequence1963},
Golomb sequence \cite{Golomb1993polyphase}, can probably improve
the performance of the algorithms.

\begin{figure}[H]
\centering{}\includegraphics[width=0.95\columnwidth]{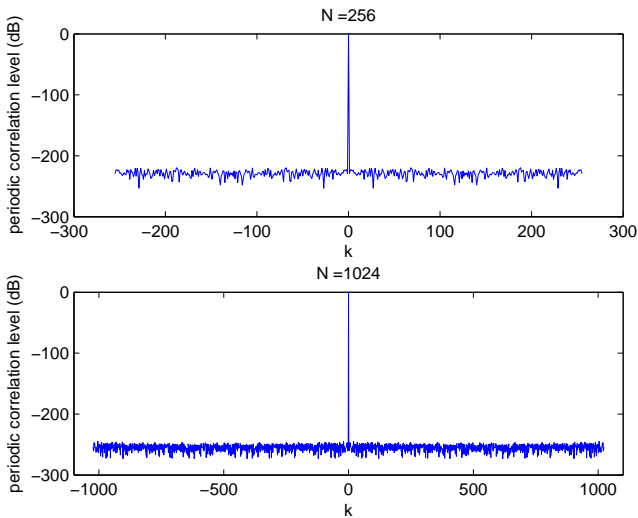}
\protect\caption{\label{fig:MISL_correlation_level-pe}Periodic correlation level of
two sequences of length $N=256$ and $N=1024$.}
\end{figure}

Finally, we present an example of applying the accelerated-MISL algorithm
to design unimodular sequences with impulse-like periodic autocorrelations.
Recall that in the periodic case, unimodular sequences with exact
impulsive autocorrelation (i.e., zero ISL) exist for any length $N$.
The periodic correlation levels of two generated sequences of length
$N=256$ and $N=1024$ are shown in Fig. \ref{fig:MISL_correlation_level-pe},
where the periodic correlation level is defined as in the aperiodic
case. 
From Fig. \ref{fig:MISL_correlation_level-pe}, we can see that the
periodic correlation levels of the two sequences are indeed very close
to zero (lower than -200dB) at lags $k\neq0$.

\subsection{Spectral-MISL}

\begin{figure}[tbh]
\centering{}\includegraphics[width=0.95\columnwidth]{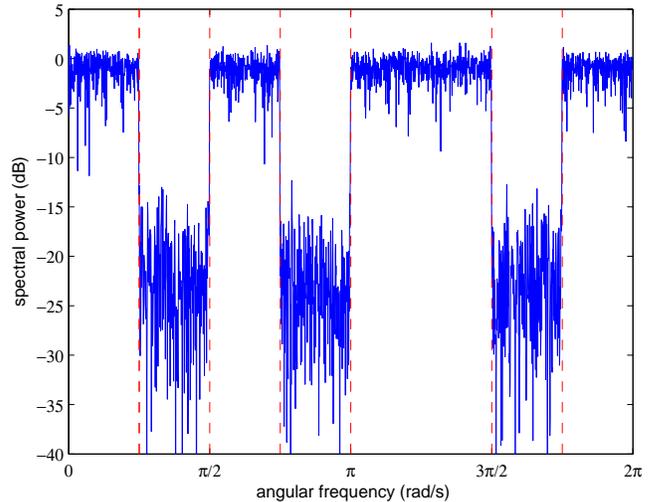}\protect\caption{\label{fig:Spectral_MISL_power}Spectral power of the designed sequence
of length $N=1000,$ generated by spectral-MISL with $\lambda=10^{4}$.}
\end{figure}

\begin{figure}[tbh]
\centering{}\includegraphics[width=0.95\columnwidth]{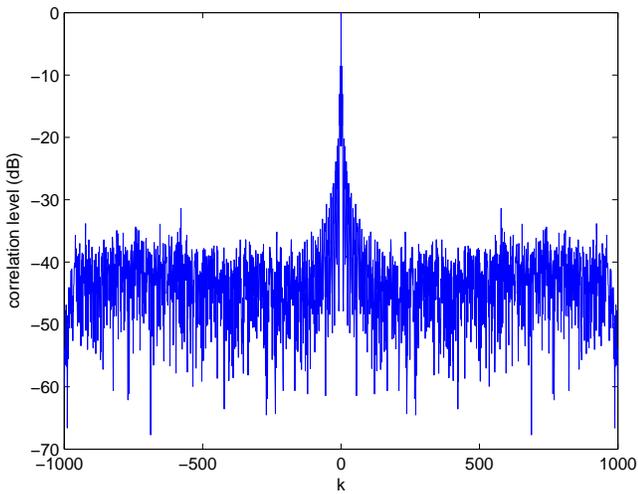}\protect\caption{\label{fig:Spectral_MISL_correlation_level}Correlation level of the
designed sequence of length $N=1000,$ generated by spectral-MISL
with $\lambda=10^{4}$.}
\end{figure}

In this subsection, we consider the design of a sequence of length
$N=1000$ with low autocorrelation sidelobes and at the same time
with low spectral power in the frequency bands $[\frac{\pi}{4},\frac{\pi}{2})\cup[\frac{3\pi}{4},\pi)\cup[\frac{3\pi}{2},\frac{7\pi}{4})$.
To achieve this, we set the index set in spectral-MISL as $\Omega=\{250,\ldots,499\}\cup\{750,\ldots,999\}\cup\{1500,\ldots,1749\}$
accordingly and tune the parameter $\lambda$ to adjust the tradeoff
between minimizing the correlation level and restricting the power
in the pre-specified frequency bands. Fig. \ref{fig:Spectral_MISL_power}
shows the spectral power of the output sequence generated by spectral-MISL
when initialized with a random sequence and with parameter $\lambda=10^{4}$.
The correlation level of the same output sequence is shown in Fig.
\ref{fig:Spectral_MISL_correlation_level}. We can see from the figures
that the power in the pre-specified frequency bands has been suppressed
and the correlation sidelobes are still relatively low.

\section{Conclusion}

We have presented an efficient algorithm named MISL for the minimization
of the ISL metric of unimodular sequences. The MISL algorithm is derived
based on the general majorization-minimization scheme and the convergence
to a stationary point is guaranteed. Acceleration schemes that can
be used to speed up MISL have also been considered. Numerical results
show that in the case of aperiodic autocorrelations the proposed MISL
algorithm can generate sequences with larger merit factor than the
state-of-the-art method and in the case of periodic autocorrelations,
MISL can produce sequences with virtually zero autocorrelation sidelobes.
A MISL-like algorithm, i.e., spectral-MISL, has also been developed
to minimize the ISL metric when there are additional spectral constraints.
Numerical experiments show that the spectral-MISL algorithm can be
used to design sequences with spectral power suppressed in arbitrary
frequency bands and at the same time with low autocorrelation sidelobes.
All the proposed algorithms can be implemented by means of the FFT
and thus are computationally very efficient in practice.

\appendix

\subsection{Proof of $\lambda_{{\rm max}}(\mathbf{\Phi})=2N^{2}$ \label{sub:Proof-Lambda_max}}
\begin{IEEEproof}
We first recall that
\begin{equation}
\mathbf{\Phi}=\sum_{p=1}^{2N}{\rm vec}(\mathbf{A}_{p}){\rm vec}(\mathbf{A}_{p})^{H},
\end{equation}
where $\mathbf{A}_{p}=\mathbf{a}_{p}\mathbf{a}_{p}^{H}$, $\mathbf{a}_{p}=[1,e^{j\omega_{p}},\cdots,e^{j\omega_{p}(N-1)}]^{T}$
and $\omega_{p}=\frac{2\pi}{2N}(p-1),\: p=1,\cdots,2N.$ The $(m,n)$-th
element of matrix $\mathbf{A}_{p}$ is given by $\left[\mathbf{A}_{p}\right]_{m,n}=e^{j\omega_{p}(m-n)},\, m=1,\ldots,N,\, n=1,\ldots,N,$
which is also the $\left(m+(n-1)N\right)$-th element of ${\rm vec}(\mathbf{A}_{p})$.
Then the $\left(m_{1}+(n_{1}-1)N,m_{2}+(n_{2}-1)N\right)$-th element
of the matrix ${\rm vec}(\mathbf{A}_{p}){\rm vec}(\mathbf{A}_{p})^{H}$
is given by $\left[{\rm vec}(\mathbf{A}_{p}){\rm vec}(\mathbf{A}_{p})^{H}\right]_{m_{1}+(n_{1}-1)N,m_{2}+(n_{2}-1)N}=e^{j\omega_{p}(m_{1}-n_{1}-m_{2}+n_{2})},\, m_{1}=1,\ldots,N,\, n_{1}=1,\ldots,N,\, m_{2}=1,\ldots,N,\, n_{2}=1,\ldots,N.$
Then it is easy to see that
\begin{align}
\left[\Phi\right]_{m_{1}+(n_{1}-1)N,m_{2}+(n_{2}-1)N} & =\sum_{p=1}^{2N}e^{j\omega_{p}(m_{1}-n_{1}-m_{2}+n_{2})}\nonumber \\
 & =\begin{cases}
2N, & m_{1}-m_{2}=n_{1}-n_{2}\\
0, & {\rm otherwise}.
\end{cases}
\end{align}
 We can see that $\mathbf{\Phi}$ is a real matrix.

To prove $\lambda_{{\rm max}}(\mathbf{\Phi})=2N^{2},$ we will show
that $\mathbf{x}^{T}\left(2N^{2}\mathbf{I}-\mathbf{\Phi}\right)\mathbf{x}\geq0,\,\forall\mathbf{x}\in\mathbf{R}^{N^{2}}$
and the equality is achievable for some $\mathbf{x}\neq\mathbf{0}.$
Let us define $\mathcal{I}_{0}=\{1,N+2,\ldots,N^{2}-N-1,N^{2}\},$
$\mathcal{I}_{-1}=\{2,N+3,\ldots,(N-1)N\},$ $\mathcal{I}_{1}=\{N+1,2N+2,\ldots,N^{2}-1\},\ldots,\mathcal{I}_{1-N}=\{N\}$
and $\mathcal{I}_{N-1}=\{N^{2}-N+1\}.$ It can be shown that 
\begin{equation}
\begin{aligned} & \mathbf{x}^{T}\left(2N^{2}\mathbf{I}-\mathbf{\Phi}\right)\mathbf{x}\\
= & 2N\left(\sum_{k=1-N}^{N-1}\sum_{i,j\in\mathcal{I}_{k},i\neq j}(x_{i}-x_{j})^{2}+\sum_{k=1-N}^{N-1}\sum_{i\in\mathcal{I}_{k}}\left|k\right|x_{i}^{2}\right).
\end{aligned}
\label{eq:xPhix}
\end{equation}
From \eqref{eq:xPhix}, it is easy to see that $\mathbf{x}^{T}\left(2N^{2}\mathbf{I}-\mathbf{\Phi}\right)\mathbf{x}\geq0,\,\forall\mathbf{x}\in\mathbf{R}^{N^{2}}$
and the equality holds for all $\mathbf{x}\in\mathbf{R}^{N^{2}}$
with the following structure: 
\[
\begin{cases}
x_{i}=x_{j}, & i,j\in\mathcal{I}_{0}\\
x_{i}=0, & {\rm otherwise}.
\end{cases}
\]
The proof is complete.
\end{IEEEproof}
\bibliographystyle{IEEEtran}
\bibliography{CAN}

\begin{thebibliography}{10}
\providecommand{\url}[1]{#1}
\csname url@samestyle\endcsname
\providecommand{\newblock}{\relax}
\providecommand{\bibinfo}[2]{#2}
\providecommand{\BIBentrySTDinterwordspacing}{\spaceskip=0pt\relax}
\providecommand{\BIBentryALTinterwordstretchfactor}{4}
\providecommand{\BIBentryALTinterwordspacing}{\spaceskip=\fontdimen2\font plus
\BIBentryALTinterwordstretchfactor\fontdimen3\font minus
  \fontdimen4\font\relax}
\providecommand{\BIBforeignlanguage}[2]{{%
\expandafter\ifx\csname l@#1\endcsname\relax
\typeout{** WARNING: IEEEtran.bst: No hyphenation pattern has been}%
\typeout{** loaded for the language `#1'. Using the pattern for}%
\typeout{** the default language instead.}%
\else
\language=\csname l@#1\endcsname
\fi
#2}}
\providecommand{\BIBdecl}{\relax}
\BIBdecl

\bibitem{turyn1968sequences}
R.~Turyn, ``Sequences with small correlation,'' \emph{Error correcting codes},
  pp. 195--228, 1968.

\bibitem{golomb2005signal}
S.~W. Golomb and G.~Gong, \emph{{Signal Design for Good Correlation: For
  Wireless Communication, Cryptography, and Radar}}.\hskip 1em plus 0.5em minus
  0.4em\relax Cambridge University Press, 2005.

\bibitem{levanon2004radar}
N.~Levanon and E.~Mozeson, \emph{{Radar Signals}}.\hskip 1em plus 0.5em minus
  0.4em\relax John Wiley \& Sons, 2004.

\bibitem{he2012waveform}
H.~He, J.~Li, and P.~Stoica, \emph{{Waveform Design for Active Sensing Systems:
  A Computational Approach}}.\hskip 1em plus 0.5em minus 0.4em\relax Cambridge
  University Press, 2012.

\bibitem{MeritFactor1972}
M.~Golay, ``A class of finite binary sequences with alternate auto-correlation
  values equal to zero (corresp.),'' \emph{IEEE Transactions on Information
  Theory}, vol.~18, no.~3, pp. 449--450, May 1972.

\bibitem{barker1953group}
R.~Barker, ``Group synchronizing of binary digital systems,''
  \emph{Communication theory}, pp. 273--287, 1953.

\bibitem{Binary2golay1977}
M.~J. Golay, ``Sieves for low autocorrelation binary sequences,'' \emph{IEEE
  Transactions on Information Theory}, vol.~23, no.~1, pp. 43--51, 1977.

\bibitem{Binary4mertens1996}
S.~Mertens, ``Exhaustive search for low-autocorrelation binary sequences,''
  \emph{Journal of Physics A: Mathematical and General}, vol.~29, no.~18, pp.
  473--481, 1996.

\bibitem{Binary1Kocabas2003}
S.~Kocabas and A.~Atalar, ``Binary sequences with low aperiodic autocorrelation
  for synchronization purposes,'' \emph{IEEE Communications Letters}, vol.~7,
  no.~1, pp. 36--38, Jan. 2003.

\bibitem{Binary5jedwab2005survey}
J.~Jedwab, ``A survey of the merit factor problem for binary sequences,'' in
  \emph{Sequences and Their Applications-SETA 2004}.\hskip 1em plus 0.5em minus
  0.4em\relax Springer, 2005, pp. 30--55.

\bibitem{Binary3wang2008}
S.~Wang, ``Efficient heuristic method of search for binary sequences with good
  aperiodic autocorrelations,'' \emph{Electronics Letters}, vol.~44, no.~12,
  pp. 731--732, 2008.

\bibitem{FrankSequence1963}
R.~L. Frank, ``Polyphase codes with good nonperiodic correlation properties,''
  \emph{IEEE Transactions on Information Theory}, vol.~9, no.~1, pp. 43--45,
  Jan. 1963.

\bibitem{Golomb1993polyphase}
N.~Zhang and S.~W. Golomb, ``Polyphase sequence with low autocorrelations,''
  \emph{IEEE Transactions on Information Theory}, vol.~39, no.~3, pp.
  1085--1089, 1993.

\bibitem{borwein2005polyphase}
P.~Borwein and R.~Ferguson, ``Polyphase sequences with low autocorrelation,''
  \emph{IEEE Transactions on Information Theory}, vol.~51, no.~4, pp.
  1564--1567, 2005.

\bibitem{Polycode2009}
C.~Nunn and G.~Coxson, ``Polyphase pulse compression codes with optimal peak
  and integrated sidelobes,'' \emph{IEEE Transactions on Aerospace and
  Electronic Systems}, vol.~45, no.~2, pp. 775--781, Apr. 2009.

\bibitem{De_Maio2009}
A.~De~Maio, S.~De~Nicola, Y.~Huang, Z.-Q. Luo, and S.~Zhang, ``Design of phase
  codes for radar performance optimization with a similarity constraint,''
  \emph{IEEE Transactions on Signal Processing}, vol.~57, no.~2, pp. 610--621,
  Feb. 2009.

\bibitem{stoica2009new}
P.~Stoica, H.~He, and J.~Li, ``New algorithms for designing unimodular
  sequences with good correlation properties,'' \emph{IEEE Transactions on
  Signal Processing}, vol.~57, no.~4, pp. 1415--1425, 2009.

\bibitem{ITROX2012}
M.~Soltanalian and P.~Stoica, ``Computational design of sequences with good
  correlation properties,'' \emph{IEEE Transactions on Signal Processing},
  vol.~60, no.~5, pp. 2180--2193, May 2012.

\bibitem{PeCAN}
P.~Stoica, H.~He, and J.~Li, ``On designing sequences with impulse-like
  periodic correlation,'' \emph{IEEE Signal Processing Letters}, vol.~16,
  no.~8, pp. 703--706, Aug. 2009.

\bibitem{he2010waveform}
H.~He, P.~Stoica, and J.~Li, ``Waveform design with stopband and correlation
  constraints for cognitive radar,'' in \emph{2nd International Workshop on
  Cognitive Information Processing (CIP)}, Elba Island, Italy, Jun. 2010, pp.
  344--349.

\bibitem{hunter2004MMtutorial}
D.~R. Hunter and K.~Lange, ``A tutorial on {MM} algorithms,'' \emph{The
  American Statistician}, vol.~58, no.~1, pp. 30--37, 2004.

\bibitem{razaviyayn2013unified}
M.~Razaviyayn, M.~Hong, and Z.-Q. Luo, ``A unified convergence analysis of
  block successive minimization methods for nonsmooth optimization,''
  \emph{SIAM Journal on Optimization}, vol.~23, no.~2, pp. 1126--1153, 2013.

\bibitem{Aldo2013}
G.~Scutari, F.~Facchinei, P.~Song, D.~P. Palomar, and J.-S. Pang,
  ``Decomposition by partial linearization: Parallel optimization of
  multi-agent systems,'' \emph{IEEE Transactions on Signal Processing},
  vol.~62, no.~3, pp. 641--656, Feb. 2014.

\bibitem{Bertsekas2003}
D.~P. Bertsekas, A.~Nedi\'{c}, and A.~E. Ozdaglar, \emph{{Convex Analysis and
  Optimization}}.\hskip 1em plus 0.5em minus 0.4em\relax Athena Scientific,
  2003.

\bibitem{SQUAREM}
R.~Varadhan and C.~Roland, ``Simple and globally convergent methods for
  accelerating the convergence of any {EM} algorithm,'' \emph{Scandinavian
  Journal of Statistics}, vol.~35, no.~2, pp. 335--353, 2008.

\bibitem{raydan2002CBB}
M.~Raydan and B.~F. Svaiter, ``Relaxed steepest descent and
  {Cauchy-Barzilai-Borwein} method,'' \emph{Computational Optimization and
  Applications}, vol.~21, no.~2, pp. 155--167, 2002.

\bibitem{barzilai1988_BB}
J.~Barzilai and J.~M. Borwein, ``Two-point step size gradient methods,''
  \emph{IMA Journal of Numerical Analysis}, vol.~8, no.~1, pp. 141--148, 1988.

\bibitem{cognitiveRadar_Haykin2006}
S.~Haykin, ``Cognitive radar: a way of the future,'' \emph{IEEE Signal
  Processing Magazine}, vol.~23, no.~1, pp. 30--40, Jan. 2006.

\end{thebibliography}

\end{document}